\documentclass[12pt,a4paper]{article}
\usepackage[english]{babel}
\usepackage[utf8]{inputenc}
\usepackage{lmodern}
\usepackage[T1]{fontenc}
\usepackage{graphicx}
\usepackage{amsmath,amsthm,amssymb}
\usepackage{amsfonts}
\usepackage{amssymb}
\usepackage{enumitem}
\usepackage{esint}
\usepackage{soul}
\usepackage{floatrow}
\usepackage[a4paper,hcentering,vcentering]{geometry}
\usepackage{subeqnarray}
\usepackage[all]{xy}
\usepackage{varioref}
\usepackage{tikz}
\usetikzlibrary{decorations.pathreplacing,calligraphy}
\usepackage{hyperref}
\usepackage{pgfplots}
\DeclareMathOperator{\cat}{\text{CAT}}
\DeclareMathOperator{\Mod}{\text{Mod}}
\DeclareMathOperator{\Out}{\text{Out}}

\DeclareMathOperator{\supp}{\text{supp}}

\DeclareMathOperator{\iso}{Isom}

\newcommand{\bd}{\partial_\infty}
\newcommand{\bdg}{\partial_{\text{Grom}}}
\newcommand{\nui}{\check{\nu}}
\newcommand{\mui}{\check{\mu}}

\newcommand{\prob}{\text{Prob}}
\theoremstyle{plain}
\newtheorem{thm}{Theorem}[section]
\newtheorem*{Rank thm}{Rank Rigidity Theorem}
\newtheorem*{Rank cjt}{Rank Rigidity Conjecture}
\newtheorem{cor}[thm]{Corollary}
\newtheorem{lem}[thm]{Lemma}
\newtheorem{prop}[thm]{Proposition}
\theoremstyle{definition} 
\newtheorem{Def}[thm]{Definition}
\newtheorem{rem}[thm]{Remark}
\theoremstyle{definition} 

\theoremstyle{definition}
\newtheorem*{ackn}{Acknowledgement}

\title{Central limit theorem on $\cat$(0) spaces with contracting isometries}
\author{Corentin Le Bars}
\date{}

\begin{document}
	\maketitle
	
	\begin{abstract}
		Let $G$ be a group with a non-elementary action on a proper $\cat(0) $ space $X$, and let $\mu$ be a measure on $G$ such that the random walk $(Z_n)_n$ generated by $\mu$ has finite second moment on $X$. Let $o$ be a basepoint in $X$, and assume that there exists a rank one isometry in $G$. We prove that in this context, $(Z_n o )_n$ satisfies a central limit theorem, namely that the random variables $\frac{1}{\sqrt{n}}(d(Z_n o, o) - n \lambda) $ converge in law to a Gaussian distribution $N_\mu$, for $\lambda$ the (positive) drift of the random walk. The strategy relies on the use of hyperbolic models introduced by H. Petyt, A. Zalloum and D. Spriano in \cite{petyt_spriano_zalloum22}, which are analogues of curve graphs and cubical hyperplanes for the class of $\cat$(0) spaces. As a side result, we prove that the probability that the nth-step $Z_n$ acts on $X$ as a contracting isometry goes to 1 as $n$ goes to infinity. 
	\end{abstract}

	\tableofcontents
	\section{Introduction}	
	
	Let $G$ be a discrete group acting by isometries on a proper $\cat$(0) space $X$. Let $\mu $ be a probability measure on $G$, which we always assume admissible, meaning that the support of $\mu$ generates $G$ as a semigroup. Consider the sequence $\omega = (\omega_i)_i$, where the $\omega_i's $ are chosen independently according to the measure $\mu$. The random walk $(Z_n(\omega))_n$ on $G$ generated by $\mu $ is then defined by $Z_n(\omega) = \omega_1 \dots \omega_n$. Taking $o \in X $, we want to study the asymptotic behaviour of the random variables $(Z_n(\omega)o )_n$. To be more precise, we want to study limit laws of the random walk in a natural compactification of $X$. Even though these questions may be hard to solve for general metric spaces, the theory is very rich when $X$ possesses nice linear or hyperbolic-like properties. In the fundamental paper of V. Kaimanovich \cite{kaimanovich00}, the convergence of $(Z_n o)_o$ to a point of the visual boundary is proven for groups acting geometrically on proper hyperbolic spaces and several other classes of actions. More recently this result has been extended by J. Maher and G. Tiozzo in \cite{maher_tiozzo18} for groups acting by isometries on non proper hyperbolic spaces. A major difficulty in the proof of the latter result was that in the non proper setting, the completion of a hyperbolic space by its Gromov boundary might be non compact. The results of Maher and Tiozzo will be fundamental in the sequel because we will deal with hyperbolic spaces without properness assumption. In \cite[Theorem 2.1]{karlsson_margulis}, Karlsson and Margulis proved a first general result of convergence of the random walk on $\cat$(0) spaces, under the assumption that the escape rate $\lambda = \liminf \frac{d(Z_no, o)}{n}$ is positive. 
	
	In \cite{LeBars22} we proved that if $G$ acts on a $\cat$(0) space $X$ with rank one isometries, then the random walk $(Z_n(\omega))_n$ almost surely converges to a point of the boundary of the visual compactification $\bd X$. A rank one element is an axial isometry whose axes do not bound any flat half plane. We give more details on this notion in Section \ref{background}, but a rank one element must be thought of as a contracting isometry with features that typically arise in hyperbolic settings. In this context, we also prove that the escape rate (the \emph{drift}) is almost surely positive: there exists $\lambda >0$ such that almost surely, $\lim_n \frac{d(Z_n o, o )}{n} = \lambda$. We review these results in Section \ref{results hyp}. The present paper can be thought of as a continuation of \cite{LeBars22}, and the goal of this work is to study further limit laws of the random walk $(Z_n o)_n$, and more specifically central limit theorems for the random variables $(d(Z_n o, o ))_n$. 
	\newline
	
	In the case of a random product of matrices, a classical result of Furstenberg \cite{furstenberg63} is the following. Take $(M_n)$ a sequence of matrices in $\text{GL}_n(\mathbb{R})$, independent and identically distributed according to a probability measure $\mu$ whose support generates a noncompact subgroup of  $\text{GL}_n(\mathbb{R})$ that does not preserve any proper linear subspace of $\mathbb{R}^n $. Assume that $\mu$ has finite first moment. Then there exists $\lambda >0$ such that for all $v \in \mathbb{R}^n - \{0\}$, 
	\begin{equation*}
		\frac{1}{n} \log \| M_n \dots M_1 v\| \longrightarrow_n \lambda
	\end{equation*}
	almost surely. This result can be thought of as an analogue of a law of large numbers on the random walk $(M_n \dots M_1 v)_n$. In this context, central limit theorems and other limit laws were proven by Furstenberg-Kesten \cite{furstenberg_kesten60}, Le Page \cite{lepage82} and Guivarc'h-Raugi \cite{guivarch_raugi85}. These state that there exists $\sigma_\mu >0 $ such that for every $v \in \mathbb{R}^n \setminus \{0\}$, 
	\begin{equation*}
		\frac{\log \| M_n \dots M_1 v\| - n  \lambda}{\sqrt{n}} \underset{n}{\longrightarrow} \mathcal{N}(0, \sigma_\mu^2), 
	\end{equation*}
	where $\mathcal{N}(0, \sigma^2)$ is a centred Gaussian law on $\mathbb{R}$. We recall that the convergence in law means that for any bounded continuous function $F : \mathbb{R} \rightarrow \mathbb{R}$, one has 
	
	\begin{equation*}
		\lim_{n \rightarrow \infty} \int_G F(\frac{\log \| M_n \dots M_1 v\| - n  \lambda}{\sqrt{n}})d\mu^{\ast n}(g) = \int_\mathbb{R} F(t) \frac{\exp(-t^2 / 2 \sigma^2)}{\sqrt{2\pi \sigma^2}}dt.
	\end{equation*}
	
	Those kinds of results were also obtained in negative curvature settings, for example in Gromov-hyperbolic groups \cite{bjorklund09}. However, the results stated thus far were proven under rather strong moment conditions. Typically, $\mu$ was assumed to have a finite exponential moment, that is, for which there exists $\alpha >0 $ such that $\int_G \exp(\alpha d(o, g o)) d\mu(g)<\infty$. 
	
	Recently, Benoist and Quint have developed a new approach to this question and have proven central limit theorems in the linear context \cite{benoist_quint16CLTlineargroups} and for hyperbolic groups \cite{benoist_quint16}. They could weaken the moment condition and only assume that the measure $\mu $ has finite second moment $\int_G (\log\|g v\|)^2 d\mu(g) < \infty$. Namely, if $\mu$ is such a measure on a group $G$ acting non elementarily on a proper hyperbolic space $Y$ with basepoint $o$, then there exists $\lambda >0$ such that the random variables $\frac{1}{\sqrt{n}}(d(Z_n(\omega)o, o) - n  \lambda)  $ converge in law to a non-degenerate Gaussian distribution \cite[Theorem 1.1]{benoist_quint16}.
	\newline

	Using this approach, C. Horbez proved central limit theorems for mapping class groups of closed connected orientable hyperbolic surfaces and on $\Out(F_N)$ \cite{horbez18}. More recently, T. Fern\'os, J. Lécureux and F. Mathéus proved that if $G$ is a group acting non-elementarily on a finite-dimensional $\cat$(0) cube complex, then we also have a central limit theorem for the random variables $(d(Z_n(\omega) o,o))_n$ \cite{fernos_lecureux_matheus21}. In both cases, the authors only assume a second moment condition. 
	
	The main result of this paper is to prove a similar result in the context of a group acting on a general $\cat$(0) space, under the assumption that an element of the group acts as a rank one isometry. We say that the group action $G \curvearrowright X$ is non-elementary if there are no fixed points in $\overline{X}$ nor a fixed pair of points in $\bd{X}$. 
	
	\begin{thm}\label{thm clt cat0}
		Let $G$ be a discrete group and $G \curvearrowright X$ a non-elementary action by isometries on a proper $\cat$(0) space $X$. Let $\mu \in \prob(G) $ be an admissible probability measure on $G$ with finite second moment, and assume that $G $ contains a rank one element. Let $o \in X$ be a basepoint of the random walk. Let $\lambda$ be the (positive) drift of the random walk. Then the random variables $\frac{1}{\sqrt{n}}(d(Z_n o, o) - n \lambda) $ converge in law to a non-degenerate Gaussian distribution $N_\mu$. 
	\end{thm}
	
	Our strategy relies heavily on the approach developed by Benoist and Quint. To summarize, one needs to approximate the random walk by a well-chosen cocycle. Then, they give a general criterium (Theorem \ref{criterium clt benoist quint} below) under which this cocycle converges in law to a Gaussian distribution. 
	
	To apply this strategy, one needs to obtain good estimates on this cocycle. The general idea of this paper is then the following. In order to get a precise description of the random walk, we use a hyperbolic space that is conveniently attached to the original $\cat$(0) space. As the theory of random walks in hyperbolic spaces is rich, we study the behavior of $\{Z_n o\}_n$ on this model, and then we lift this information back to the original $\cat$(0) space. This strategy was implemented successfully in \cite{horbez18} and \cite{fernos_lecureux_matheus21}:
	
	\begin{itemize}
		\item for $\Mod(S)$, the hyperbolic model is the curve complex $C(S)$, and the lifting to $\mathcal{T}(S)$ is done in \cite[Section 3.4]{horbez18}; 
		\item for a $\cat$(0) cube complex, the hyperbolic model is the contact graph $\mathcal{C}X$, and the lifting is implemented in \cite[section 5]{fernos_lecureux_matheus21}. 
	\end{itemize}
	
	In \cite{petyt_spriano_zalloum22}, H. Petyt, D. Spriano and A. Zallum introduced analogues of curve graphs and cubical hyperplanes for the class of $\cat$(0) spaces. Using a generalized notion of hyperplane, they build a family of hyperbolic metrics $(d_L)_L$ on $X$ which conserve many of the geometric features of the original $\cat$(0) space. These spaces capture hyperbolic behaviours in $X$ and behave very well under the isometric action of a group. Moreover, a rank one isometry of $X$ acts on some hyperbolic model as a loxodromic isometry. Our strategy will be to chose a good hyperbolic model $X_L = (X, d_L)$, and then to make use of the limit laws proven by Maher and Tiozzo in \cite{maher_tiozzo18}, or by Gou\"ezel in \cite{gouezel22}. A key fact is that there is an equivariant homeomorphic embedding of the Gromov boundary $\bdg X_L$ of the hyperbolic model $X_L$ into the visual boundary of the $\cat$(0) space \cite[Theorem 7.1]{petyt_spriano_zalloum22}. 
	\newline
	
	Another interesting question in the study of $(Z_n(\omega))_n$ is the proportion of steps that are "hyperbolic". In the context of random walks on hyperbolic spaces, Maher and Tiozzo show that the probability that a random walk of size $n$ is a loxodromic isometry goes to $1$ as $n$ goes to infinity \cite[Theorem 1.4]{maher_tiozzo18}. For a non-elementary action on an irreducible $\cat$(0) cube complex, Fern\'os, Lécureux and Mathéus show that the proportion of steps $Z_n$ that are contracting goes to $1$ as $n$ goes to infinity. They use this result to show that if a group $G$ acts non-elementarily and essentially on a (possibly reducible) finite-dimensional $\cat$(0) cube complex, then there exist regular elements, extending a result of Caprace and Sageev \cite{caprace_sageev11}. In our context, we also prove that "most" of the steps in the random walk are rank one. This result is not involved in the proof of Theorem \ref{thm clt cat0}, but is of independent interest. 
	
	\begin{thm}[Rank one elements in the random walk]
		Let $G$ be a discrete group and $G \curvearrowright X$ a non-elementary action by isometries on a proper $\cat$(0) space $X$. Let $\mu \in \prob(G) $ be an admissible probability measure on $G$, and assume that $G $ contains a rank one element. Then 
		\begin{eqnarray}
			\mathbb{P}(\omega \, : \, Z_n (\omega) \text{ is a contracting isometry }) \underset{n\rightarrow \infty}{\rightarrow} 1 \nonumber.
		\end{eqnarray}
	\end{thm}
	Using the curtain models from \cite{petyt_spriano_zalloum22}, such a result is actually straightforward. We emphasize the idea that a systematic approach of dynamics on $\cat$(0) spaces using these hyperbolic models can prove fruitful, especially when quantitative estimates are required. Indeed, curtain models also benefit from their combinatorial structure. In this paper, we exploit this richness in several ways, especially in the main geometric lemma \ref{geometric lemma}. In \cite[Section 5]{le-bars23these}, the author develops these connections in order to study other limit laws on general Hadamard spaces. 
	\newline 
	
	A different approach for the study of such limit laws was implemented in \cite{mathieu_sisto20}, where the authors prove central limit theorems on acylindrically hyperbolic groups. Their strategy relies on a control of deviation inequalities, which encapsulate the way the random walk progresses in an ``almost aligned'' way, hence their approach apply to possibly non-proper spaces. While there is a slight overlap with the results stated here (especially \cite[Theorem 13.4]{mathieu_sisto20}), Mathieu and Sisto study random walks on acylindrically hyperbolic groups with a word metric. This situation does not immediately apply here. Indeed, in our main theorem, the pull-back metric induced on G by an orbit map need not be quasi-isometric to a word metric, and in fact need not even be proper. Also, their assumptions on the measure $\mu$ are much more restrictive: they assume that $\mu$ has finite exponential moment. In particular, their assumption is not optimal, while it is the case here. Last, the techniques involved are completely different: their approach has a ``local'' flavor, whereas here we use boundary theory and compactifications. 
	\newline

	While we were working on this project, Inhyeok Choi released a paper in which he states central limit theorems along with other limit laws in $\cat$(0) spaces, Teichmüller spaces and outer spaces \cite{choi22a}. One of the main assumptions is still the presence of a pair of independent contracting isometries in the group, but the methods and the proofs are different. Indeed, Choi uses a pivotal technique introduced by Boulanger, Mathieu and Sisto in \cite{mathieu_sisto20} and \cite{boulanger_mathieu_sisto21} and further developed by Gouëzel in \cite{gouezel22}. These techniques have a ``local flavor'', while our paper relies on boundary theory, and uses hyperbolic models that depend on specific features of $\cat$(0) spaces. We think this approach is natural from a geometric point of view, and we believe that the interplay between $\cat$(0) spaces and their underlying hyperbolic models will be useful in the study of still open questions about limit laws. In any case, it is always interesting to have different strategies and techniques for studying radom walks and limit laws. 
	\newline
	
	The essential assumption in these results is the presence of contracting elements for the action of $G$ on the $\cat$(0) space $X$. In \cite[Chapter 5]{le-bars23these}, the author proves that actually, all the results presented here hold in the more general context of a Hadamard space, i.e. a separable and complete $\cat$(0) space, removing the properness assumption on $X$. Notice that the boundary of a general Hadamard spaces may be no longer compact, and that the embedding of the boundaries $\bdg X_L \hookrightarrow \bd X$ is only stated for proper $\cat$(0) spaces in \cite{petyt_spriano_zalloum22}. In \cite[Theorem 5.3.5]{le-bars23these}, we prove that this embedding actually holds in this more general setting. Once this is done, the general strategy for proving the central limit theorem \ref{thm clt cat0} is similar, although with some additional technical difficulties. We refer to the manuscript \cite[Chapter 5]{le-bars23these} for details. 
	\newline 
	
	We believe our approach can be of use in order to determine if the boundary $\bd X$ endowed with the hitting measure is actually the Poisson boundary of $(G, \mu)$, extending a result of Karlsson and Margulis for cocompact actions \cite[Corollary 6.2]{karlsson_margulis}. 
	
	Moreover, it seems natural to use these hyperbolic spaces to prove that if $\mu$ has finite first moment, then limit points of the random walk almost surely belong to the sublinear Morse boundary constructed by Qing and Rafi in \cite{qing_rafi22}. Note that this question is linked to the previous one, because it is believed that the sublinear Morse boundary is often a good candidate for the Poisson boundary, especially for finitely supported measures, see for example \cite[Theorem F]{qing_rafi22} and \cite[Theorem B]{qing_rafi_tiozzo20}. In both cases, the use of hyperbolic models seems useful because of precision of estimates that can be derived from the combinatorial structure of these spaces. 
	\newline
	
	In Section \ref{background}, we review basic definitions about random walks, rank one isometries and explain our setting. In Section \ref{hyperbolic models}, we 
	explicit the construction and properties of the hyperbolic models $(X, d_L)$, and give various geometric lemmas that will be useful afterwards. Section \ref{results hyp} is dedicated to presenting the works of Maher and Tiozzo in \cite{maher_tiozzo18} and of Gou\"ezel in \cite{gouezel22}, and the first results in proper $\cat(0) $ spaces that were found in \cite{LeBars22}. We explain the strategy developed by Benoist and Quint in Section \ref{strategy}, and give the proof of our main Theorem in Section \ref{Section proof}. 
	
	\begin{ackn}
		This work was done during the trimester "Groups acting on fractals, hyperbolicity and self-similarity", which was held at IHP from April 11th to July 8th, 2022. The author acknowledges support of the Institut Henri Poincaré (UAR 839 CNRS-Sorbonne Université), and LabEx CARMIN (ANR-10-LABX-59-01). We thank Jean Lécureux for his help and commentaries on this paper. The author is also thankful to A. Sisto for drawing our attention on another approach for these questions, and to I. Choi for discussions on this matter. We are very grateful to H. Petyt, D. Spriano and A. Zalloum for allowing us to discuss the details of their construction and for their friendly remarks. We are grateful to the anonymous referee for pointing out some mistakes and giving advice on how to improve the presentation of this paper. 
	\end{ackn}
	
	\section{Background}\label{background}
	
	\subsection{Random walks and $\cat$(0) spaces}
	
	Let $G$ be a discrete countable group and $\mu \in \prob(G)$ a probability measure on $G$. Recall that the support of $\mu$ is 
	\begin{equation*}
		\supp(\mu) := \{g \in G \, | \, \mu (g ) > 0\}. 
	\end{equation*}
	\begin{Def}
		We say that a measure $\mu$ on a discrete countable group is \emph{admissible} if its support $\supp(\mu)$ generates $G$ as a semigroup.
	\end{Def}
	Throughout the article we will assume that $\mu $ is admissible. Let $(\Omega, \mathbb{P}) $ be the probability space $(G^{\mathbb{N}}, \delta_e \times \mu^{\mathbb{N^\ast}})$, where $\delta_e$ is the Dirac measure at $e$. The application 
	\begin{equation*}
		(n, \omega) \in \mathbb{N} \times \Omega \mapsto Z_n(\omega) = \omega_1 \omega_2 \dots \omega_n,
	\end{equation*}
	where $\omega$ is chosen according to the law $\mathbb{P}$, defines the random walk on $G$ generated by the measure $\mu$.

	Let now $(X,d)$ be a proper $\cat$(0) metric space, on which $G$ acts by isometries. If the reader wants a detailed introduction to $\cat (0) $ spaces, the main references that we will use are \cite{bridson_haefliger99} and \cite{ballman95}. We recall that the boundary $ \bd X$ of a $\cat (0) $ space $X$ is the set of equivalent classes of rays $\sigma : [0, \infty) \rightarrow X$, where two rays $\sigma_1, \sigma_2$ are equivalent if they are asymptotic, i.e. if $d(\sigma_1(t), \sigma_2 (t))$ is bounded uniformly in $t$. 
	
	Given two points on the boundary $\xi$ and $\eta$, if there exists a geodesic line $\sigma : \mathbb{R} \rightarrow X$ such that the geodesic ray $\sigma_{[0, \infty)}$ is in the class of $\xi$ and the geodesic ray $t \in [0, \infty) \mapsto \sigma(-t)$ is in the class of $\eta$, we will say that the points $\xi $ and $\eta$ are joined by a geodesic line. The reader should be aware that in general, such a geodesic need not exist between any two points of the boundary, as can be seen in $\mathbb{R}^2$. A point $\xi$ of the boundary is called a \textit{visibility point} if, for all $\eta \in \bd X - \{\xi\}$, there exists a geodesic from $\xi$ to $\eta$. We will see in the next section a criterion to prove that a given boundary point is a visibility point. 
	
	An important feature in $\cat$(0) spaces is the existence of closest-point projections on complete convex subsets. More precisely, given a complete convex subset $C$ in a $\cat$(0) space, there exists a map $\pi_C : X \rightarrow C$ such that $\pi_C(x)$ minimizes the distance $d(x,C)$: 
	
	\begin{prop}[{\cite[Lemma 2.4]{bridson_haefliger99}}]\label{projection prop}
		The projection $\pi_C$ onto a convex complete subset in a $\cat$(0) space satisfies the following properties:
		\begin{itemize}
			\item $\forall x \in X$, $\pi_C(x)$ is uniquely defined and $d(x, \pi_C(x))= d(x, C) = \inf_{c \in C} d(x, c) $; 
			\item if $x'$ belongs to the geodesic segment $[x, \pi(x)]$, then $\pi_C(x')= \pi_C(x)$; 
			\item $\pi_C $ is a retraction of $X$ onto $C$ that does not increase the distances: for all $x, y \in X$, we have $d(\pi_C(x), \pi_C(y )) \leq d(x,y)$.
		\end{itemize}
	\end{prop}
	
	It is immediate to see that the above properties can be applied to geodesic segments, which are convex and complete with the induced metric. When $\gamma : [a, b] \rightarrow X$ is a geodesic segment, we will write $\pi_\gamma$ for the projection onto the image $[\gamma(a), \gamma(b)]\subseteq X$.

	When $X$ is a proper space, the space $\overline{X} = X \cup \partial X $ is a compactification of $X$, that is, $\overline{X}$ is compact and $X$ is an open and dense subset of $\overline{X}$. We recall that the action of $G$ on $X$ extends to an action on $\bd X$ by homeomorphisms. 
	
	Another equivalent construction of the boundary can be done using horofunctions. If $x_n \rightarrow \xi \in \bd X$ and $x \in X$, we denote by $b_\xi^{x} : X \mapsto \mathbb{R}$ the horofunction given by 
	\begin{equation*}
		b_\xi^{x}(z) = \lim_n d(x_n, z) - d(x_n, x).
	\end{equation*}
	It is a standard result in $\cat$(0) geometry (see for example \cite[Proposition II.2.5]{ballman95}) that this limit exists and that given any basepoint $x$, a horofunction characterizes the boundary point $\xi $. When the context is clear we will often omit the basepoint and just write $ b_\xi $.
	
	\subsection{Rank one elements}\label{rank one section}
	
	Let $g \in G$. We say that $g$ is a \textit{semisimple} isometry if its displacement function $ x \in X \mapsto \tau_g(x) = d(x , gx)$ has a minimum in $X$. If this minimum is non-zero, it is a standard result (see for example \cite[Proposition II.3.3]{ballman95}) that the set on which this minimum is obtained is of the form $C \times \mathbb{R}$, where $C$ is a closed convex subset of X. On the set $\{c\}\times \mathbb{R}$ for $c \in C$, $g$ acts as a translation, which is why $ g$ is called \textit{axial} and the subset $\{c\}\times \mathbb{R}$ is called an \textit{axis} of $g$. A \textit{flat half-plane} in $X$ is defined as a euclidean half plane isometrically embedded in $X$.
	
	\begin{Def}
		We say that a geodesic in $X$ is \textit{rank one} if it does not bound a flat half-plane. If $g$ is an axial isometry of $X$, we say that $g$ is rank one if no axis of $g$ bounds a flat half-plane. 
	\end{Def}

	If $G$ acts on $X$ by isometries and possesses a rank one element $g \in G$ for this action, we may say that $G$ is rank one. However, the theory of $\cat$(0) groups is not as clear as for Gromov hyperbolic groups. For example, there is no good (i.e. invariant under quasi isometry) notion of boundary of a $\cat$(0) group, as shown by Croke and Kleiner in \cite{croke_kleiner00}. To summarize, it is better to keep in mind that "rank one" is always attached to a given action $G \curvearrowright X$ on a $\cat$(0) space. 
	\newline
	
	More information on rank one isometries and geodesics can be found in \cite[Section III. 3]{ballman95}, and more recently in \cite{caprace_fujiwara10} and in \cite{bestvina_fujiwara09}. 
	
	\begin{Def}
		We say that the action $G \curvearrowright X$ of a rank one group $G$ on a $\cat (0) $ space $X$ is \textit{non-elementary} if $G$ neither fixes a point in $\bd X$ nor stabilizes a geodesic line in $X$. 
	\end{Def}
	
	To justify this definition, we use a result from Caprace and Fujiwara in \cite{caprace_fujiwara10}. What follows comes from the aforementioned paper. 
	
	\begin{Def}
		Let $g_1, \, g_2 \in G$ be axial isometries of $G$, and fix $x_0 \in X$. The elements $g_1, g_2 \in G$ are called independent if the map 
		\begin{equation}
			\mathbb{Z} \times \mathbb{Z} \rightarrow [0, \infty) : (m,n) \mapsto d(g_1^m x_0, g_2^nx_0)
		\end{equation}
		is proper. 
	\end{Def}
	
	\begin{rem}
		In particular, the fixed points of two independent axial elements form four distinct points of the visual boundary. 
	\end{rem}
	
	Let us end this section by stating two results about rank one isometries. The first one was proven by P-E.~Caprace and K.~Fujiwara in \cite{caprace_fujiwara10}. 
	
	\begin{prop}[{\cite[Proposition 3.4]{caprace_fujiwara10}}]\label{non elem caprace fuj}
		Let $X$ be a proper $\cat $(0) space and let $G < \iso (X)$. Assume that $G$ contains a rank one element. Then exactly one of the following assertions holds: 
		\begin{enumerate}
			\item \label{alt elem} $G$ either fixes a point in $\bd X$ or stabilizes a geodesic line. In both cases, it possesses a subgroup of index at most 2 of infinite Abelianization. Furthermore, if $X$ has a cocompact isometry group, then $\overline{G} < \iso (X)$ is amenable. 
			
			\item \label{alt non elem} G contains two independent rank one elements. In particular, $\overline{G}$ contains a discrete non-Abelian free subgroup. 
		\end{enumerate}
	\end{prop}
	
	As a consequence, the action $G \curvearrowright X$ of a rank one group $G$ on a $\cat (0) $ space $X$ is non-elementary if and only if alternative \ref{alt non elem} of the previous Proposition holds. 
	\newline
	
	Rank one isometries are especially interesting because they induce natural contracting properties on the space. These properties mimic how loxodromic isometries behave in the hyperbolic setting. 
	
	\begin{Def}
		A geodesic $\sigma$ in a $\cat$(0) space is said to be \textit{$C$-contracting} with $C >0$ if for every metric ball $B$ disjoint from $\sigma$, the projection $\pi_\sigma (B)$ of the ball $B$ onto $\sigma$ has diameter at most $C$. An axial isometry is contracting if there exists $C>0$ such that one of its axes is $C$-contracting. 
	\end{Def} 
	
	It is clear that a contracting isometry is rank one. It turns out that the converse is true if $X$ is a proper $\cat$(0) space, as was shown by M. Bestvina and K. Fujiwara in \cite{bestvina_fujiwara09}. This result will allow us to use the hyperbolic models described in Section \ref{hyperbolic models}. 
	
	\begin{thm}[{\cite[Theorem 5.4]{bestvina_fujiwara09}}]\label{bestv fuj rank one contracting}
		Let $X$ be a proper $\cat$(0) space, $g : X \rightarrow X$ be an axial isometry and $\sigma$ be an axis of $g$. Then there exists $B$ such that $\sigma$ is $B$-contracting if and only if $\sigma$ does not bound a half-flat. In other words, $g$ is contracting if and only if $g$ is a rank one isometry. 
	\end{thm}

	\subsection{Gromov products}
	
	Let $(X,d) $ be a metric space. One defines the Gromov product of $x, y \in X$ with respect to $o \in X$ as 
	\begin{eqnarray}
		(x |y )_o = \frac{1}{2}(d(x, o) + d(y, o) - d(x,y)). \nonumber
	\end{eqnarray}
	The quantity $(x |y )_o$ must be thought of as representing the distance between $o$ and the geodesic between $x$ and $y$. This notion is particularly interesting because it does not require $X$ to be actually geodesic, and in fact we often deal with only quasigeodesic spaces. Also, we can use Gromov products to characterize hyperbolic spaces. We recall that a metric space $(X,d) $ is hyperbolic if there is $\delta >0$ such that for all $x, y, z \in X$, 
	
	\begin{eqnarray}
		(x|z)_o \geq \min((x|y)_o, (y|z)_o) - \delta \nonumber.
	\end{eqnarray}
	
	If the reader wants a detailed introduction to hyperbolic spaces, a standard reference is \cite{bridson_haefliger99}.
	\newline
	
	If $(X,d)$ is a proper $\cat(0) $ space, the Gromov product can be extended to the visual boundary $\bd X$ of $X$ by the following formulas: for $x, y \in \bd X$, $o, m \in X$, 
	\begin{eqnarray}
		(m|x)_o & := &\frac{1}{2}(d(o, m) - b^o_x(m)); \nonumber \\
		(x |y)_o & := &- \frac{1}{2}\, \underset{q \in X}{\inf}\, (b_x (q) + b_{y}(q)) \label{eq gromov product}.
	\end{eqnarray}
	Assume that $x , y \in X$. A quick computation shows that the infimum of Equation~\eqref{eq gromov product} is attained for any $q\in [x,y]$. Indeed, for any other $p \in X$, 
	\begin{eqnarray}
		b_x (p) + b_{y}(p) - (b_x (q) + b_{y}(q)) & = & d(x, p) + d(y, p) - (d(o, x) + d(o , y)) \nonumber \\
		& &- (d(x, q ) + d(y, q)) + (d(o, x) + d(o , y)) \nonumber \\
		& = & d(x, p) + d(y, p) - d(x, y) \text{ because $q \in [x,y]$} \nonumber \\
		& \geq & 0 \text{ by the triangular inequality}. \nonumber
	\end{eqnarray}
	When $x, y  \in \bd X$, take $(x_n) $ and $(y_n)$ sequences converging to $x,y \in \bd X $ respectively for the visual topology. Since the visual compactification is equivalent to the compactification by Busemann functions (see for instance \cite[Theorem II.8.13]{bridson_haefliger99}), $\{b_{x_n}\}$ and $\{b_{y_n}\} $ converge to $b_x$ and $b_y$ respectively (for the topology of uniform converge on bounded sets). In particular, if there exists a geodesic line $\gamma$ such that $\gamma(t ) \underset{t \to \infty}{\to} x$ and $\gamma(t ) \underset{t \to -\infty}{\to} y$, then the previous computation shows that for any point $q \in \gamma$:
	\begin{eqnarray}
		(x |y)_o & = &  \lim_{n, m \to \infty} - \frac{1}{2}\,  (b_{\gamma(n)} (q) + b_{\gamma(-m)}(q)) \nonumber \\
		& = & \lim_{n,m \to \infty} (x_n |y_m)_o \nonumber.
	\end{eqnarray} 
	
	\section{Hyperbolic models for proper $\cat$(0) spaces}\label{hyperbolic models}

	The goal of this section is to briefly present some ideas of \cite{petyt_spriano_zalloum22}, in which the authors build a way of attaching a family of hyperbolic metric spaces $X_L = (X, d_L)_L$ to a proper $\cat$(0) space. What is interesting about these spaces is that they convey much of the geometry of the original space, especially at infinity, and they behave very well under isometric actions. More specifically, rank one isometries will act on some well-chosen spaces as loxodromic isometries. This construction can be understood as the analogue (and generalization) of the curve graphs that exist in the context of $\cat$(0) cube complexes, see \cite{hagen14} and \cite{genevois19}.

	\begin{Def}
		Let $X$ be a $\cat(0) $ space, and let $\gamma: I \rightarrow X$ be a geodesic. Let $\pi_\gamma$ be the projection onto the geodesic $\gamma$ characterized by Proposition \ref{projection prop}. Let $t \in I$ be such that $[t - \frac{1}{2} , t + \frac{1}{2}]$ belongs to $I$. Then the \textit{curtain} dual to $\gamma$ at $t$ is 
		\begin{eqnarray}
			h = h_{\gamma, t} = \pi^{-1}_\gamma (\gamma ([t - \frac{1}{2} , t + \frac{1}{2}])). \nonumber
		\end{eqnarray}
		The \textit{pole} of $h_{\gamma, t}$ is $ \gamma ([t - \frac{1}{2} , t + \frac{1}{2}]) $. Borrowing from the vocabulary of hyperplanes, we will call $h^{-} =  \pi^{-1}_\gamma (\gamma ((-\infty, t- \frac{1}{2})\cap I))$ and $h^{+} =  \pi^{-1}_\gamma (\gamma ((t +\frac{1}{2}, + \infty )\cap I))$ the \textit{halfspaces} determined by $h$. Note that $\{h^{-} , h , h^{+}\}$ is a partition of $X$. If $A \subseteq h^{-}$ and $B \subseteq h^{+}$ are subsets of $X$, we say that $h$ \textit{separates} $A$ from $B$. 	
	\end{Def}

	We will often denote a curtain by the letter $h$, even though one must keep in mind that $h = h_{\gamma, t}$ is characterized by a given geodesic $\gamma : I \rightarrow X $ and a point $t \in I$ (which defines a unique pole $P \subseteq \gamma$). Sometimes, we may also write $h= h_{\gamma, P}$ to emphasize on the pole $P$. 
		
	\begin{rem}\label{remark curtains thick}
		By Proposition \ref{projection prop}, it is immediate that curtains are closed subsets of $X$, and that they are thick: if $h$ is a curtain, then $d(h^{-}, h^{+})= 1$. 
	\end{rem}
	
	Curtains can fail to be convex: if $x, y \in h^{-}$, it may happen that there exists $z \in [x, y ]\cap h^{+}$, see \cite[Remark 3.4]{petyt_spriano_zalloum22}. Nonetheless, we have a weaker notion of convexity that the authors call star convexity:
	
	\begin{prop}[{\cite[Lemma 2.6]{petyt_spriano_zalloum22}}]\label{star convexity}
		Let $h$ be a curtain dual to $\gamma$ and $P\subseteq \gamma$ be its pole. For every $x \in h$, then $[x, \pi_P (x)] \subseteq h$. 
	\end{prop}
	
	\begin{Def}
		A family of curtains $\{h_i\}$ is said to be a chain if $h_i $ separates $h_{i-1}$ from $h_{i+1}$ for every $i $. Chains can be used in order to define a metric on $X$ by the following: for $x \neq y \in X$, 
		\begin{eqnarray}
			d_\infty (x, y) = 1 + \max \{\, |c| \, : \, c \text{ is a chain separating }x \text{ from } y \}. \nonumber
		\end{eqnarray}
	\end{Def}
	
	One can check that this definition gives a metric. If $h$ is a curtain, we have seen that $d(h^{-}, h^{+}) = 1$, hence for any $x,y \in X$, $d_\infty (x, y ) \leq \lceil d(x,y ) \rceil$. Conversely, it turns out that $d$ and $d_\infty$ may differ by at most 1, as shown by the following lemma. 
	
	\begin{lem}[{\cite[Lemma 2.10]{petyt_spriano_zalloum22}}]
		Let $x, y \in X$. Then there is a chain of curtains $c$ dual to $[x,y]$ that realizes $d_\infty (x, y) = 1 + |c| $. and for which $ 1 + |c| = \lceil d(x,y) \rceil$. 
	\end{lem}
	
	We are now ready to refine the notion of separation in order to capture only some of the hyperbolic features of the space. 
	
	We say that a chain $c$ of curtains \emph{meets} a curtain $h$ if every single curtain $h_i \in c$ intersects $h$. 
	
	\begin{Def}[$L$-separation]
		Let $L \in \mathbb{N}^\ast$, we say that disjoint curtains are $L$\textit{-separated} if every chain meeting both has cardinality at most $L$. A chain of pairwise $L$-separated curtains is called an $L$\textit{-chain}. 
	\end{Def}
	
\begin{figure}
		\centering
		\begin{center}
			\begin{tikzpicture}[scale=1.2]
				\draw (0,0) -- (4,0)  ;
				\draw (2.8,2) -- (2.8,-1)  ;
				\draw (1.2,2) -- (1.2,-1)  ;
				\draw (2,2) -- (2,-1)  ;
				\draw (2,0.5) node[right]{$\leq 2L+1$} ;
				\draw (-1, 2.5) to[bend right = 80] (5, 3);
				\draw (0,0) node[below left]{$x_1$} ;
				\draw (2, -1) node[below]{$h_2$} ;
				\draw (1.2, -1) node[below]{$h_1$} ;
				\draw (2.8, -1) node[below]{$h_3$} ;
				\draw (-1, 2.5) node[above left]{$x_2$} ;
				\draw (4,0) node[below right]{$y_1$} ;	
				\draw (5,3) node[above]{$y_2$};
				\filldraw[black] (2,0) circle(1pt);
				\filldraw[black] (2,1) circle(1pt);
			\end{tikzpicture}
		\end{center}
		\caption{Illustration of Lemma \ref{bottleneck}.}\label{bottleneck figure}
	\end{figure}
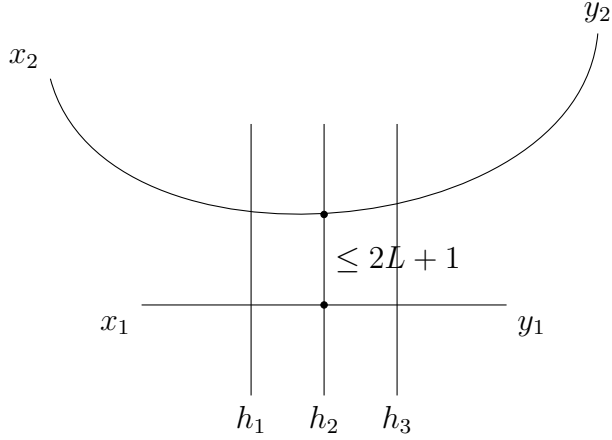

	The following geometric Lemma is a key ingredient for the proof of Theorem \ref{theorem hyperbolic models}, and will be used several times in the sequel. It means that $L$-separation induces good Morse properties. The picture one has to keep in mind is given by Figure \ref{bottleneck figure}. 
	
	\begin{lem}[{\cite[Lemma 2.14]{petyt_spriano_zalloum22}}]\label{bottleneck}
		Suppose that $A$, $B$ are two sets which are separated by an L-chain $\{h_1, h_2, h_3\}$ all of whose elements are dual to a geodesic $\gamma= [x_1, y_1]$ with $x_1 \in A$ and $y_1 \in B$. Then for any $x_2 \in A$, $y_2 \in B$, if $p \in h_2 \cap [x_2, y_2]$, then $d(p, \pi_\gamma(p))\leq 2L + 1$. 
	\end{lem}
	
	The next Lemma states that if there is a $L$-chain separating two points $x $ and $y$, we can find a (smaller) $L$-chain of curtains separating those, which is dual to the geodesic $[x,y]$ and whose size can be controlled. It will prove useful later on, especially when we want to use Lemma \ref{bottleneck}. 
	
	\begin{lem}[{\cite[Lemma 2.21]{petyt_spriano_zalloum22}}]\label{dual chain}
		Let $L, n \in \mathbb{N}$, and let $\{h_1, \dots, h_{(4L +10)n}\}$ be an $L$-chain separating $A$, $B \subseteq X$. Take $x \in A$, $y \in B$. Then $A$ and $B$ are separated by an $L$-chain of size $\geq n+1 $ dual to $[x,y]$. 
	\end{lem}
	
	We are now ready to define a family of metrics using $L$-separation. 
	
	\begin{Def}
		Given distinct points $x\neq y \in X$, we define 
		\begin{eqnarray}
			d_L(x,y) = 1 + \max\{|c| \, : \, c \text{ is an } L\text{-chain separating }x \text{ from } y\}. \nonumber
		\end{eqnarray}
	\end{Def}
	
	It turns out that for every $L$, $d_L$ gives a metric on $X$ \cite[Lemma 2.17]{petyt_spriano_zalloum22}. We will denote by $X_L = (X, d_L)$ the resulting metric space. With this definition in hand, Petyt, Spriano and Zalloum prove that the metric spaces $(X, d_L)$ are hyperbolic. 
	
	\begin{thm}[{\cite[Theorem 3.1]{petyt_spriano_zalloum22}}] \label{theorem hyperbolic models}
		For any $\cat$(0) space $X$ and any integer $L$, the space $(X, d_L)$ is a quasi-geodesic hyperbolic space with hyperbolicity constants depending only on $L$. Moreover, $\iso(X)$ acts by isometries on $(X, d_L)$. 
	\end{thm}
	
	We will then call $(X, d_L)$ a hyperbolic model for the $\cat$(0) space $X$. Another useful fact about these spaces is that they behave well under isometries with "hyperbolic-like" properties.
	
	\begin{thm}[{\cite[Theorem 4.9]{petyt_spriano_zalloum22}}] \label{rank one loxodromic}
		Let $g$ be a semisimple isometry of $X$. The following are equivalent: 
		\begin{enumerate}
			\item $g$ is a contracting isometry of the $\cat$(0) space $X$;
			\item there exists $L \in \mathbb{N}$ such that $g$ acts loxodromically on $X_L$.
		\end{enumerate} 
	\end{thm} 

	Another piece of information brought by this construction is the relation between the Gromov boundaries $\partial X_L$ of the hyperbolic models $X_L = (X, d_L)$ and the visual boundary of the original $\cat$(0) space $(X, d)$. 
	
	\begin{Def}
		We say that a geodesic ray $\gamma : [0, \infty) \rightarrow X$ \textit{crosses} a curtain $h $ if there exists $t_0 \in [0, \infty)$ such that $h$ separates $\gamma(0) $ from $\gamma ([t_0, \infty))$. Alternatively, we may say that $h$ separates $\gamma(0)$ from $\gamma(\infty)$. Similarly, we say that a geodesic line $\gamma: \mathbb{R} \rightarrow X$ crosses a curtain $h$ if there exist $t_1, t_2 \in \mathbb{R}$ such that $h$ separates $\gamma ((-\infty, t_1])$ from $\gamma ([t_2, \infty))$. We say that $\gamma $ crosses a chain $c = \{h_i\}$ if it crosses each individual curtain $h_i$. 
	\end{Def}
	
	As a consequence of Lemma \ref{bottleneck} and Lemma \ref{dual chain}, if two geodesic rays with the same starting point cross an infinite $L$-chain $c$, then they are asymptotic, and hence equal. 
	
	\begin{rem}\label{infinite dual chain}
		Since curtains are not convex, it is not obvious that any geodesic ray $\gamma$ meeting a given curtain $h$ must cross it ($\gamma$ could meet $h$ infinitely often). However, by \cite[Corollary 3.2]{petyt_spriano_zalloum22} if $\gamma$ is a geodesic ray that meets every element of an infinite $L$-chain $c = \{h_i\}_{i \in \mathbb{N}}$, then $\gamma$ must cross $c$: for every $i$, there exists $t_i \in [0, \infty)$ such that $h_i$ separates $\gamma(0) $ from $\gamma ([t_i, \infty))$.
	\end{rem}
	
	Given $o \in X$, we define $\mathcal{B}_L $ as the subspace of $\bd X$ consisting of all geodesic rays $\gamma : [0, \infty) \rightarrow X$ starting from $o$ and such that there exists an infinite $L$-chain crossed by $\gamma$. In the case of the contact graph associated to a $\cat$(0) cube complex $X$, we had the existence of an $\iso(X)$-equivariant embedding of the boundary of the contact graph into the Roller boundary $\partial_{\mathcal{R}}X$. The following result is the analogue in the context of $\cat$(0) spaces. 
	
	\begin{thm}[{\cite[Theorem 7.1]{petyt_spriano_zalloum22}}] \label{thm equivariant embedding of boundaries}
		Let $X$ be a proper $\cat$(0) space. Then, for every $L \in \mathbb{N}^\ast$, the identity map $\iota : X \longrightarrow X_L$ induces an $\iso(X)$-equivariant homeomorphism $\partial_L : \mathcal{B}_L \longrightarrow \partial X_L$. 
	\end{thm}
	
	Recall that the \textit{support} of a Borel measure $m$ on a topological space $Y$ is the smallest closed set $C$ such that $m(Y \setminus C)= 0$. In other words $y \in \supp(m)$ if and only if for all $U $ open containing $y $, $m(U) >0$.
	
	\begin{Def}
		We say that the action by isometries of a group $G$ on a hyperbolic space $Y$ (not assumed to be proper) is \textit{non-elementary} if there are two loxodromic isometries with disjoint fixed points on the Gromov boundary. A probability measure $\mu$ on $G$ is said to be \textit{non-elementary} if its support generates a group acting non-elementarily on $Y$. 
	\end{Def}
	
	In order to use the results concerning random walks in hyperbolic spaces, we must show that the action of a group $G$ on a proper $\cat$(0) space with rank one isometries induces a non-elementary action on some hyperbolic model $(X, d_L)$. 
	
	\begin{prop}\label{prop non elem}
		Let $G$ be a group acting non-elementarily by isometries on a proper $\cat$(0) space $(X,d)$, and assume that $G$ possesses a rank one element for this action. Then there exists $L \in \mathbb{N}$ such that $G$ acts on the hyperbolic space $(X, d_L)$ non-elementarily by isometries. 
	\end{prop}

	\begin{proof}
		The action $G \curvearrowright (X,d)$ is non elementary and contains a rank one element, hence by Theorem \ref{non elem caprace fuj} there exist two independent rank one isometries $ g,h$ in $G$. By Theorem \ref{bestv fuj rank one contracting}, those rank one isometries are $B$-contracting for some $B$. Now, applying Theorem \ref{rank one loxodromic}, there exists $L \in \mathbb{N}$ such that $g$ and $h$ act on $(X, d_L)$ as loxodromic isometries. As $g $ and $h$ are independent, their fixed points form four distinct points of the visual boundary $\bd X$. Now seen in $X_L = (X, d_L)$, their fixed points sets must also form four distinct points of $\partial X_L$ because of the homeomorphism $\partial_L : \mathcal{B}_L \longrightarrow \partial X_L$. This means that the action $G \curvearrowright X_L$ is non-elementary. 
	\end{proof}

\section{Random walks and hyperbolicity}\label{results hyp}
	
	The results of Section \ref{hyperbolic models} allow us to read some information about the random walk in the hyperbolic models $X_L = (X, d_L)$, and then translate this information back to the original $\cat$(0) space. As the theory of random walks on hyperbolic spaces is well-studied, one may hope that this process is fruitful.

	\subsection{Random walks on hyperbolic spaces}

	In this section, we summarize what is known concerning random walks in hyperbolic spaces. Most of the work for the non-proper case was done by Maher and Tiozzo in \cite{maher_tiozzo18}. The first result is the convergence of the random walk to the Gromov boundary. 
	
	\begin{thm}[{\cite[Theorem 1.1]{maher_tiozzo18}}] \label{mt conv}
		Let $G$ be a countable group of isometries of a separable hyperbolic space $Y$. Let $\mu$ be a non-elementary probability distribution on $G$, and $o \in Y$ a basepoint. Then the random walk $(Z_n (\omega)o )_n $ induced by $\mu$ converges to a point $z^{+} (\omega)\in \bd X$, and the resulting hitting measure is the unique $\mu$-stationary measure on $\bd X$. 
	\end{thm}
	\begin{rem}\label{rem gouezel cv non-sep}
		Note that the previous result is stated for separable hyperbolic spaces, while in our case, the hyperbolic models are not separable. However, Gou\"ezel shows in \cite[Theorem 1.3]{gouezel22} that this result of convergence remains true for possibly non-separable hyperbolic spaces. 
	\end{rem}
	
	Assume that the measure $\mu $ has finite first moment $\int d(go, o) d\mu(g) < \infty$. Let us define the drift (or escape rate) of the random walk.  
	
	\begin{Def}
		The \textit{drift} of the random walk $(Z_n  o)_n$ on a metric space $(Y, d)$ is defined as 
		\begin{eqnarray}
		l(\mu): = \inf_n \int_\Omega d(Z_n(\omega) o, o)d\mathbb{P}(\omega) = \inf_n \int_G d(g o, o)d\mu^{\ast n}(g) \nonumber
		\end{eqnarray}
		if $\mu $ has finite first moment, and $l(\mu) := \infty$ otherwise. 
	\end{Def}
	
	If $\mu $ has finite first moment, then a classical application of Kingmann subadditive Theorem sows that
	\begin{eqnarray}
		l(\mu) = \lim_n \frac{1}{n}d(Z_n (\omega)o , o), \nonumber
	\end{eqnarray}
	and the above limit is essentially constant and finite.

	In the context of a group acting on a hyperbolic space, Gou\"ezel proves that the drift is almost surely positive with no moment condition. This can be seen as a law of large numbers. 
	
	\begin{thm}[{\cite[Theorems 1.1 and 1.2]{gouezel22}}]\label{mt drift}
		Let $G$ be a countable group of isometries of a hyperbolic space $(Y,d_Y)$. Let $\mu$ be a non-elementary probability distribution on $G$, and $o \in Y$ a basepoint. Then the drift $l(\mu):= \lim_n  \frac{1}{n}d(Z_n o, o )$ is well-defined, essentially constant and positive (possibly infinite). 
		
		Moreover, for every $r < l(\mu)$, there exists $\kappa >0$ such that 
		\begin{eqnarray}
			\mathbb{P}\big( \omega \in \Omega \, : \, d_Y(Z_n(\omega) o, o) \leq rn \big)< e^{-\kappa n }.
		\end{eqnarray}
	\end{thm}

	Another piece of information that can be given about the random walk is the proportion of hyperbolic isometries in the random variables $(Z_n)_n$. Recall that the translation length of an isometry in a hyperbolic space is defined as $|g| := \lim_n \frac{1}{n} d(g^n o , o ) $, which does not depend on the basepoint $o$. 
	
	\begin{thm}[{\cite[Theorem 1.4]{maher_tiozzo18}}]\label{mt 4}
		Let $G$ be a countable group of isometries of a separable hyperbolic space $Y$. Let $\mu$ be a non-elementary probability distribution on $G$, and $o \in Y$ is a basepoint. Then the translation length $|Z_n(\omega)| $ grows almost surely at least linearly in $ n $: there exists $K>0 $ such that
		\begin{eqnarray}
			\mathbb{P}(\omega \, : \, |Z_n(\omega)| \leq Kn) \underset{n \rightarrow \infty}{\longrightarrow} 0 \nonumber.
		\end{eqnarray} 
	\end{thm}

	The above result thus implies that the probability that $Z_n(\omega)$ is not a loxodromic isometry goes to zero as $n$ goes to infinity. 
	
	\subsection{First results for random walks in $\cat$(0) spaces}
	
	In $\cat(0) $ spaces, many of the previous theorems hold if we assume that there are elements in the acting group $G$ that share "hyperbolic-like" properties. Namely, if $X$ is a proper $\cat$(0) space, we will assume that $G$ contains rank one isometries of $X$. The first result deals with stationary measures on $\overline{X}$. Recall that a measure $\nu \in \prob(\overline{X})$ is called stationary if $\mu \ast \nu = \nu$. 
	
	\begin{thm}[{\cite[Theorem 1.1]{LeBars22}}]\label{measure thm cat0}
		Let $G$ be a discrete group and $G \curvearrowright X$ a non-elementary action by isometries on a proper $\cat $(0) space $X$. Let $\mu \in \prob(G) $ be an admissible probability measure on $G$, and assume that $G $ contains a rank one element. Then there exists a unique $\mu$-stationary measure $\nu \in \prob (\overline{X})$. 
	\end{thm}

	 The convergence of the random walk to the boundary can then be established in this setting. It is the analogue of Theorem \ref{mt conv}.
	
	\begin{thm}[{\cite[Theorem 1.2]{LeBars22}}]\label{convergence thm}
		Let $G$ be a discrete group and $G \curvearrowright X$ a non-elementary action by isometries on a proper $\cat$(0) space $X$. Let $\mu \in \prob(G) $ be an admissible probability measure on $G$, and assume that $G $ contains a rank one element. Then for every $x \in X$, and for $\mathbb{P}$-almost every $\omega \in \Omega$, the random walk $(Z_n (\omega) x)_n $ converges almost surely to a boundary point $z^{+}(\omega) \in\bd X$. Moreover, $z^{+}(\omega)$ is distributed according to the stationary measure $\nu$. 
	\end{thm} 
	
	Interestingly, we can prove that the limit points are almost surely rank one, meaning that for almost any pair of limit points $\xi, \eta \in \bd X$, there exists a rank one geodesic in $X$ joining $\xi$ to $\eta$ (\cite[Corollary 1.3]{LeBars22}). This feature suggests the use of hyperbolic models. First, we establish a result concerning the proportion of rank one elements in the random walk. 
	
	\begin{thm}\label{thm prop rank one elements}
		Let $G$ be a discrete group and $G \curvearrowright X$ a non-elementary action by isometries on a proper $\cat$(0) space $X$. Let $\mu \in \prob(G) $ be an admissible probability measure on $G$, and assume that $G $ contains a rank one element. Then 
		\begin{eqnarray}
			\mathbb{P}(\omega \, : \, Z_n (\omega) \text{ is a contracting isometry }) \underset{n\rightarrow \infty}{\rightarrow} 1 \nonumber.
		\end{eqnarray}
	\end{thm}

	\begin{proof}
		Because of Proposition \ref{prop non elem}, we can then apply the results of Maher-Tiozzo and Gou\"ezel. In particular, by Theorem \ref{mt 4}, the translation length $|Z_n(\omega)|_L$ of $(Z_n(\omega))_n$ grows almost surely at least linearly in $n$. Therefore, the probability that $Z_n(\omega)$ is a loxodromic element of $X_L$ goes to 1 as $ n $ goes to $\infty$. But thanks to Theorem \ref{rank one loxodromic}, an isometry $g$ of the $\cat$(0) space $X$ is contracting if and only if there is an $L $ such that $g$ acts as a loxodromic isometry on $X_L$. The previous argument now implies that the probability that $Z_n(\omega)$ is a contracting isometry of $X$ goes to 1 as $n $ goes to $ \infty$. 
	\end{proof}
	
	\begin{rem}
		There is a slight omission in the proof of Theorem \ref{thm prop rank one elements}. Indeed, Theorem \ref{mt 4} is stated for geodesic, separable hyperbolic spaces, while hyperbolic models are non-separable and only almost geodesic. However, thanks to Bonk and Schramm \cite[Theorem 4.1 ]{bonk_schramm00}, this result extends to non-geodesic hyperbolic spaces. The key thing is that due to Theorem \ref{mt drift} we have control of the displacement variables $d(Z_n o, o) $ up to the escape rate $l(\mu)$. A detailed proof of Theorem \ref{thm prop rank one elements} can be found in \cite[Section 5.3.5]{le-bars23these}.
	\end{rem}
	
	The analogue of Theorem \ref{mt drift} also holds in the context of $\cat$(0) spaces with rank one isometries. 
	
	\begin{thm}[{\cite[Theorem 1.4]{LeBars22}\label{drift thm}}]
		Let $G$ be a discrete group and $G \curvearrowright X$ a non-elementary action by isometries on a proper $\cat$(0) space $X$. Let $\mu \in \prob(G) $ be an admissible probability measure on $G$ with finite first moment, and assume that $G $ contains a rank one element. Let $o \in X$ be a basepoint of the random walk. Then the drift $\lambda$ is almost surely positive: 
		\begin{equation}
			\lim_{n \rightarrow \infty} \frac{1}{n} d(Z_n o, o) = \lambda >0. \nonumber
		\end{equation}
	\end{thm}

	Actually H. Izeki worked on the drift-free case in \cite{izeki22}. The author proves a strengthening of Theorem \ref{drift thm}, in that it is valid even for finite dimensional, non proper $\cat$(0) spaces, and without the assumption that there are rank one isometries. The counterpart is that one needs to assume that $\mu$ has finite second moment. Namely, Izeki proves that in this context, either the drift $\lambda$ is strictly positive, or there is a $G$-invariant flat subspace in $X$ \cite[Theorem A]{izeki22}. However, for our purpose, we will only need Theorem \ref{drift thm}. 
	\newline
	
	In the proof of Theorem \ref{drift thm}, we actually show that the displacement $d(Z_n(\omega) x, x)$ is almost surely well approximated by the Busemann functions $b_\xi (Z_n(\omega)x)$. This result will be used later when we give geometric estimates for the action. 
	
	\begin{prop}[{\cite[Proposition 5.2]{LeBars22}}]\label{approx displacement horo}
		Let $G$ be a discrete group and $G \curvearrowright X$ a non-elementary action by isometries on a proper $\cat $(0) space $X$. Let $\mu \in \prob(G) $ be an admissible probability measure on $G$ with finite first moment, and assume that $G $ contains a rank one element. Let $x \in X$ be a basepoint. Then for $\nu$-almost every $\xi  \in \partial X$, and $\mathbb{P}$-almost every $\omega \in \Omega$, there exists $C > 0 $ such that for all $n \geq 0$ we have 
		\begin{equation}
			|b_\xi (Z_n(\omega)x) - d(Z_n(\omega) x, x)| < C. 
		\end{equation}
	\end{prop}

	\section{Central Limit Theorems and general strategy} \label{strategy}
	
	In order to prove our main result, we use a strategy that is largely inspired by the works of Benoist and Quint on linear spaces and hyperbolic spaces, see \cite{benoist_quint16} and \cite{benoist_quint16CLTlineargroups}. They developed a method for proving central limit theorems for cocycles, relying on results due to Brown in the case of martingales \cite{brown71}. 
	
	\subsection{Centerable cocycle} 
	
	Let $G$ be a discrete group, $Z$ a compact $G$-space and $c$ a cocycle $c : G \times Z \rightarrow \mathbb{R}$, meaning that $c(g_1g_2, x) = c(g_1, g_2 x ) + c(g_2, x)$, and assume that $c$ is continuous. Let $\mu$ be a probability measure on $G$. 
	
	\begin{Def}
		Let $c$ be a continuous cocycle $c : G \times Z \rightarrow \mathbb{R}$. We say that $c$ has \textit{constant drift} $c_\mu$ if $c_\mu = \int_G c(g, x) d\mu(g)$ does not depend on $x \in Z$. We say that $c$ is \textit{centerable} if there exists a bounded measurable map $\psi : Z \rightarrow \mathbb{R} $ and a cocycle $c_{0} : G \times Z \rightarrow \mathbb{R}$ with constant drift $c_{0, \mu} = \int_G c_0(g, x) d\mu(g)$ such that 
		\begin{eqnarray}
			c(g, x) = c_0(g,x) + \psi(x) - \psi(gx).
		\end{eqnarray}
		We say that $c$ and $c_0$ are cohomologous. In this case, the \textit{average} of $c$ is defined to be $c_{0, \mu}$. 
	\end{Def}

	\begin{rem}\label{rem average cocycle}
			Let $\nu \in \prob(Z)$ be a $\mu$-stationary measure, and let $c : G \times Z \rightarrow \mathbb{R}$ be a centerable continuous cocycle: for $g \in G, x \in Z$, $c(g, x) = c_0(g,x) + \psi(x) - \psi(gx)$ with $c_0$ having constant drift and $\psi$ bounded measurable. The following computation shows that the average of $c$ does not depend on the particular choices of $c_0$ and $\psi$. Indeed: 
		\begin{eqnarray}
			\int_{G\times Z} c(g, x) d\mu(g) d\nu (x) &  = & \int_{G\times Z} c_0(g, x) d\mu(g) d\nu (x) + \int_Z \psi(x) d\nu(x)  \nonumber \\
			& & - \int_{G \times Z} \psi(gx) d\mu(g) d\nu(x) \nonumber \\
			& = & \int_{G} c_0(g, x) d\mu(g) + \int_Z \psi(x) d\nu(x)  - \int_{G \times Z} \psi(gx) d\mu(g) d\nu(x) \nonumber \\
			& = &  \int_{G} c_0(g, x) d\mu(g) + \int_Z \psi(x) d\nu(x) \nonumber \\
			& &- \int_{ Z} \psi(x) d\nu(x) \text{ because $\nu$ is $\mu$-stationary} \nonumber \\
			& = & c_{0, \mu} \text{ because $c_0$ has constant drift}. \nonumber
		\end{eqnarray}
		Hence the average of $c$ is given by $\int c(g, x) d\mu(g) d\nu(x) $, which explains the terminology. Moreover, the average of $c$ does not depend on the choices of $c_0$ and $\psi$. 
	\end{rem}
	
	The reason why we study limit laws on cocycles is the following result. This version is borrowed from Benoist and Quint, who improved previous results from Brown about central limit theorems for martingales \cite{brown71}.
	
	\begin{thm}[{\cite[Theorem 3.4]{benoist_quint16CLTlineargroups}}]\label{criterium clt benoist quint}
		Let $G$ be a locally compact group acting by homeomorphisms on a compact metrizable space $Z$. Let $c : G \times Z \rightarrow \mathbb{R}$ be a continuous cocycle such that $\int_G \sup_{x \in Z} | c(g,x)|^2 d\mu (g) < \infty $. Let $\mu $ be a Borel probability measure on $G$. Assume that $c$ is centerable with average $\lambda_c$ and that there exists a unique $\mu$-stationary probability measure $\nu $ on $Z$. 
		
		Then the random variables $\frac{1}{\sqrt{n}}(c(Z_n, x) - n \lambda_c) $ converge in law to a Gaussian law $N_\mu$. In other words, 
		for any bounded continuous function $F$ on $\mathbb{R}$, one has
		\begin{eqnarray}
			\int_G F \big( \frac{c(g,x) - n \lambda_c}{\sqrt{n}}\big) d(\mu^{\ast n})(g) \longrightarrow \int_{\mathbb{R}} F(t) dN_\mu (t) \nonumber. 
		\end{eqnarray}
		Moreover, if we write $c (g, z) = c_0(g, z) + \psi(z) - \psi(g z) $ with $\psi $ bounded and $c_0$ with constant drift $c_\mu$, then the covariance 2-tensor of the limit law is 
		\begin{eqnarray}
			\int_{G \times Z} (c_0(g, z) - c_\mu)^2 d\mu(g) d\nu(z) \nonumber. 
		\end{eqnarray}
	\end{thm}

	\subsection{Busemann cocycle and strategy}
	
	Let $G$ be a discrete group and $G \curvearrowright X$ a non-elementary action by isometries on a proper $\cat$(0) space $X$. Let $\mu \in \prob(G) $ be an admissible probability measure on $G$ with finite first moment, and assume that $G $ contains a rank one element. Let $o \in X$ be a basepoint of the random walk. Theorems \ref{convergence thm} and \ref{drift thm} ensure that the random walk $(Z_n(\omega)o)_n$ converges to a point of the boundary and that the drift $\lambda = \lim_n \frac{1}{n} d(Z_n(\omega) o , o)$ is well-defined and almost surely positive. 
	
	We denote by $\mui$ the probability measure on $G$ defined by $\mui(g) = \mu(g^{-1})$. Let $(\check{Z}_n)_n$ be the right random walk associated to $\mui$. Since $\mu $ is admissible and has finite first moment, so does $\mui$. We can then apply Theorems \ref{measure thm cat0}, \ref{convergence thm} and \ref{drift thm} to $\mui$. We will denote by $\nui $ the unique $\mui$-stationary measure on $\overline{X}$, and by $\check{\lambda} $ the positive drift of the random walk $(\check{Z}_n o)_n$. 
	
	\begin{rem}
		One can check that
		\begin{eqnarray}
			\check{\lambda} &=& \inf_n \frac{1}{n} \int d(go, o) d\mui^{\ast n}(g) \nonumber \\
			& = & \inf_n \frac{1}{n} \int d(o, g^{-1} o) d\mui^{\ast n}(g) \nonumber \\
			& = & \inf_n \frac{1}{n} \int d(o, g o) d\mu^{\ast n}(g) \nonumber,
		\end{eqnarray}
		hence $\lambda = \check{\lambda}$.
	\end{rem}
	
	In our context, the continuous cocycle that we consider is the Busemann cocycle on the visual compactification of the $\cat$(0) space $X$: for $x \in \overline{X}, \ g \in G$ and $o \in X$ a basepoint, 
	
	\begin{equation*}
		\beta(g,x) = b_x (g^{-1} o).
	\end{equation*}

	It is straightforward to show that $\beta$ is continuous. Observe that for all $g_1, g_2 \in G, \, x \in Y$, horofunctions satisfy a cocycle relation: 
	\begin{eqnarray}
		b_\xi (g_1g_2 o ) & = &  \lim_{x_n \rightarrow \xi}d(g_1g_2 , x_n) - d(x_n, x)  \nonumber \\
		& = & \lim_{x_n \rightarrow \xi}d(g_2 , g_1^{-1}x_n) - d(g_1 o, x_n ) +   d(g_1 o, x_n ) - d(x_n, o) \nonumber \\ 
		& = & \lim_{x_n \rightarrow \xi}d(g_2x , g_1^{-1}x_n) - d(o, g_1^{-1}x_n ) +   d(g_1 x, x_n ) - d(x_n, o) \nonumber \\ 
		& = & b_{g_1^{-1}\xi} (g_2 o) + b_\xi (g_1 o). \label{cocycle horof}
	\end{eqnarray}
	By \eqref{cocycle horof}, $\beta$ satisfies the cocycle relation $\beta(g_1g_2, x) = \beta(g_1, g_2 x ) + \beta(g_2, x)$. Thanks to Proposition \ref{approx displacement horo}, for every $o \in X$, for $\nu$-almost every $x  \in \partial X$, and $\mathbb{P}$-almost every $\omega \in \Omega$, there exists $C > 0 $ such that for all $n \geq 0$ we have 
	\begin{equation}\label{approx busemann displacement}
		|\beta (Z_n(\omega)^{-1}, x) - d(Z_n(\omega) o,o)| < C. 
	\end{equation}
	
	Equation \eqref{approx busemann displacement} shows that the cocycle $\beta(Z_n(\omega),x)$ "behaves" like $d(Z_n(\omega) o,o)$. Thus it makes sense to try and apply Theorem \ref{criterium clt benoist quint} to the Busemann cocycle $\beta(g,x)$. 
	\newline 
	
	Henceforth, we will assume that $\mu$ is an admissible probability measure on $G$ with finite second moment $\int_{G} d(go, o)^2 d\mu(g) < \infty$.
	
	The following proposition summarizes some properties of the Busemann cocycle. It shows that obtaining a central limit theorem on $\beta$ will imply our main result. 
	
	\begin{prop}\label{Busemann cocyle prop}
		Let $G$ be a discrete group and $G \curvearrowright X$ a non-elementary action by isometries on a proper $\cat $(0) space $X$. Let $\mu \in \prob(G) $ be an admissible probability measure on $G$ with finite second moment, and assume that $G $ contains a rank one element. Let $o \in X$ be a basepoint of the random walk. Let $\lambda$ be the (positive) drift of the random walk, and $\beta : G \times \overline{X} \rightarrow \mathbb{R}$ be the Busemann cocycle $\beta(g, x) = b_x(g^{-1} o)$. Then 
		\begin{enumerate}
			\item $\int_{G} \sup_{x \in \overline{X}}|\beta (g, x)|^2 d\mu(g) < \infty$ and $\int_{G} \sup_{x \in \overline{X}}|\beta (g, x)|^2 d\mui(g) < \infty$;\label{moment}
			\item For $\nu$-almost every $\xi \in \bd X$, $\lambda = \lim_n \frac{1}{n} \beta (Z_n (\omega), \xi)$ $\mathbb{P}$-almost surely; \label{drift cocycle}
			\item $\mathbb{P}$-almost surely, $\lambda = \int_{G \times \overline{X}} \beta(g, x) d\mu(g)d\nu(x) = \int_{G \times \overline{X}} \beta(g, x) d\mui(g)d\nui(x)$. \label{average Busemann}
		\end{enumerate}
	\end{prop}

	\begin{proof}
		As a consequence of Proposition \ref{approx displacement horo}, equation \eqref{approx busemann displacement} gives that for $\nu$-almost every $x  \in \partial X$, and $\mathbb{P}$-almost every $\omega \in \Omega$, there exists $C > 0 $ such that for all $n \geq 0$ we have 
		\begin{eqnarray}
			|\beta (Z_n(\omega)^{-1}, x) - d(Z_n(\omega) o,o)| < C. \label{eq beta distance}
		\end{eqnarray}
		Because the action is isometric and $\mu$ has finite second moment $\int_{G} d(g o, o)^2 d\mu(g)~<~\infty$, we obtain 
		\begin{eqnarray}
			\int_{G} \sup_{x \in \overline{X}}|\beta (g, x)|^2 d\mu(g)<\infty. \nonumber
		\end{eqnarray} 
		With the same argument: 
		\begin{eqnarray}
			\int_{G} \sup_{x \in \overline{X}}|\beta (g, x)|^2 d\mui(g) < \infty. \nonumber
		\end{eqnarray} 
	
		Now thanks to Theorem \ref{drift thm}, the variables $\{\frac{1}{n} d(Z_n(\omega) o , o)\}$ converge almost surely to $\lambda >0$. Since the action is isometric, we immediately get that $$\frac{1}{n} d(Z_n(\omega) o , o)~\to_n~\lambda$$ almost surely. Again, because the action is isometric, Equation \eqref{eq beta distance} tells that for $\nu$-almost every $x  \in \partial X$, and $\mathbb{P}$-almost every $\omega \in \Omega$, there exists $C > 0 $ such that for all $n \geq 0$ we have 
		\begin{eqnarray}
			|\beta (Z_n(\omega), x) - d(Z_n(\omega)^{-1} o,o)| < C.  \nonumber
		\end{eqnarray} 
		Combining these results, we obtain that for $\nu$-almost every $\xi \in \overline{X}$, and $\mathbb{P}$-almost every $\omega \in \Omega$, 
		\begin{eqnarray}
			\lambda = \lim_n \frac{1}{n} \beta (Z_n (\omega), \xi). \nonumber
		\end{eqnarray}
		
		The ideas in the proof of \ref{average Busemann} are classical. We give the details for the convenience of the reader.

		Let $T : (\Omega \times \overline{X}, \mathbb{P} \times \nui ) \rightarrow (\Omega \times \overline{X}, \mathbb{P} \times \nui )$ be defined by $T(\omega, \xi ) \mapsto (S\omega, \omega_0^{-1} \xi)$, with $S((\omega_i)_{i \in \mathbb{N}}) = (\omega_{i+1})_{i \in \mathbb{N}}$ the usual shift on $\Omega$. By \cite[Proposition 5.4]{LeBars22}, $T$ preserves the measure $\mathbb{P}\times \check{\nu}$  and is an ergodic transformation. Define $H : \Omega \times \overline{X} \rightarrow \mathbb{R}$ by 
		\begin{equation*}
			H(\omega, \xi) = h_\xi (\omega_0 o) = \beta(\omega_0^{-1}, \xi). 
		\end{equation*}
		By \ref{moment}, it is clear that $\int |H(\omega, \xi)| d\mathbb{P}(\omega)d\nui(\xi) < \infty$. 
		
		By cocycle relation \eqref{cocycle horof} one gets that 
		\begin{eqnarray}
			b_\xi (Z_n o) =  \sum_{k=1}^{n} h_{Z_k^{-1}\xi} (\omega_k o )  = \sum_{k=1}^{n} H(T^k (\omega, \xi)) \label{transient cocycle}. 
		\end{eqnarray}
	
		Then $\beta(Z_n(\omega)^{-1}, \xi) = \sum_{k=1}^{n} H(T^k (\omega, \xi))$, and by \ref{drift cocycle}, 
		\begin{eqnarray}
			\lambda = \lim_n \frac{1}{n} \sum_{k=1}^{n} H(T^k (\omega, \xi)).
		\end{eqnarray}
		
		Now, by Birkhoff ergodic theorem, one obtains that almost surely, 
		\begin{eqnarray}
			\lambda &=& \int_{\Omega \times \overline{X}} H(\omega, \xi) d\mathbb{P}(\omega)d\nui(x). \nonumber \\
			& = & \int_{\Omega \times \overline{X}} h_\xi(\omega_0 o) d\mathbb{P}(\omega)d\nui(x) \nonumber \\
			& = & \int_{G \times \overline{X}} \beta(g^{-1}, \xi) d\mu(g)d\nui(x) \nonumber \\
			& =& \int_{G \times \overline{X}} \beta(g, \xi) d\mui(g)d\nui(x)
		\end{eqnarray} 
		
		The previous computations can be done similarly for $\mu $ and $\nu $, hence we also have that 
		\begin{eqnarray}
			\lambda = \int_{G \times \overline{X}} \beta(g, x) d\mu(g)d\nu(x). \nonumber
		\end{eqnarray} 
	\end{proof}

	In order to apply Theorem \ref{criterium clt benoist quint} on the Busemann cocycle $\beta$, it remains to show that $\beta$ is centerable. If this is the case, by \ref{average Busemann} and Remark \ref{rem average cocycle}, its average must be the positive drift $\lambda$. In other words, we need to show that there exists a bounded measurable function $\psi : \overline X \rightarrow \mathbb{R}$ such that the cocycle
	\begin{eqnarray}
		\beta_0(g, x) = \beta_(g,x) - \psi(x) + \psi(gx) \nonumber
	\end{eqnarray} 
	has constant drift, so that the cohomological equation
	\begin{eqnarray}
		\beta(g, x) = \beta_0(g,x) + \psi(x) - \psi(gx). \label{cohom equation}
	\end{eqnarray}
	is verified. Then, proving the Central Limit Theorem in our context amounts to finding such a $\psi $ that is well defined and bounded. This will be done by using a hyperbolic model that can give nice estimates on the random walk.

	\section{Proof of the Central Limit Theorem}\label{Section proof}

	\subsection{Geometric estimates}

	In this section, we prove our main Theorem, following the strategy explained in Section \ref{strategy}. First, we will provide geometric estimates on the random walk that will be used later on. This is where we use the specific contraction properties provided by the curtains and the hyperbolic models discussed in Section \ref{hyperbolic models}. The goal is ultimately to prove that the candidate $\psi$ for the cohomological equation is bounded. 
	
	Let $G$ be a discrete group and $G \curvearrowright X$ a non-elementary action by isometries on a proper $\cat$(0) space $X$, and assume that $G $ contains a rank one element. Let $o \in X$ be a basepoint of the random walk. Recall that $B_L$ is defined to be the subspace of $\bd X$ consisting of all geodesic rays $\gamma : [0, \infty) \rightarrow X$ starting from $o$ and such that there exists an infinite $L$-chain crossed by $\gamma$. By Theorem \ref{thm equivariant embedding of boundaries}, there exists an $\iso(X)$-equivariant embedding $\mathcal{I} : \partial X_L \rightarrow \bd X$, whose image lies in $\mathcal{B}_L$. 
	
	\begin{prop}\label{geometric lemma}
		Let $(g_n)$ be a sequence of isometries of $G$, and let $o \in X$, $x, y \in \bd X$. Assume that there exists $\lambda, \varepsilon, A >0 $ such that:
		\begin{enumerate}[label= (\roman*)]
			\item $ \{g_n o \}_n $ converges in $(\overline{X_L}, d_L)$ to a point of the boundary $z_L \in \partial X_L$, whose image in $\bd X $ by the embedding $\mathcal{I}$ is not $y$; \label{assum cv}
			\item $d_L(g_n o, o) \geq An$; \label{assum drift hyp}
			\item $|b_x (g_n^{-1} o ) - n \lambda | \leq \varepsilon n$; \label{assum gnx}
			\item $|b_y (g_n o ) - n \lambda | \leq \varepsilon n$; \label{assum y}
			\item $|d(g_n o , o) - n \lambda | \leq \varepsilon n$. \label{assum drift cat}
		\end{enumerate}	
		Then, one obtains: 
		\begin{enumerate}
			\item $(g_nx |g_n o )_o \geq (\lambda - \varepsilon) n$; \label{geom estimate 1}
			\item $(y | g_n o )_o \leq \varepsilon n$.\label{geom estimate 2} 
		\end{enumerate}
		If moreover $A \geq 2(4L+10) \varepsilon$, then we have: 
		\begin{enumerate}[resume]
			\item $(y| g_n x )_o \leq \varepsilon n + (2L+ 1)$. \label{geom estimate 3}
		\end{enumerate}
	\end{prop}
	Before getting into the proof of this proposition, let us give an idea of what it represents. Assumptions \ref{assum cv} and \ref{assum drift hyp} express that the sequence $ \{g_n o \}_n $ converge to the visual boundary following a ``contracting direction'', with control on the size of the $L$-chain that separates $o$ from $g_n o$. Assumptions \ref{assum gnx}, \ref{assum y} are to be seen as "the distance between $y$ and $g_n o $ grows linearly" and "the distance between $x$ and $g_n^{-1} o $ grows linearly" (even though $x$ and $y$ are boundary points). Assumption \ref{assum drift cat} simply means that the average escape rate of $\{g_n o \}$ is close to $\lambda$. Proposition \ref{geometric lemma} states that in these circumstances, we can control the quantities $(g_nx |g_n o )_o$, $(y | g_n o )_o$ and $(y| g_n x )_o$, represented by the size of the dashed segments $E_1$, $E_2$ and $E_3$ respectively in Figure \ref{figure geom interpretation estimate}.  In Proposition \ref{prop centerable cocycle}, we shall need in particular the geometric estimate \ref{geom estimate 3} in order to prove that the candidate $\psi $ for the cohomological equation \eqref{cohom equation} is bounded, see the parallel with Proposition \ref{estimates} below.
	\begin{figure}[h!]
		\centering
		\begin{center}
			\begin{tikzpicture}[scale=1]
				\draw (0,0) -- (8,0)  ;
				\draw (8,0)  [-stealth] to(10,0)  ;
				\draw (2.5,4) -- (2.5,-0.5)  ;
				\draw (3,4) -- (3,-0.5)  ;
				\draw (3.5,4) -- (3.5,-0.5)  ;
				\draw (4.5,4) -- (4.5,-0.5)  ;
				\draw (5,4) -- (5,-0.5)  ;
				\draw (4,4) -- (4,-0.5)  ;
				\draw (5.5, 4) -- (5.5,-0.5)  ;
				\draw (0,0) node[below left]{$o$} ;
				\draw (9,4) .. controls (1,1) ..(-1, 4);
				\draw (-1, 4) .. controls (1,1) .. (8,0);
				\draw (8,0) .. controls (6.5,1.5) ..(9,3.9);
				\draw [decorate,
				decoration = {brace}] (2.4,4.1) --  (5.6, 4.1);
				\draw (0,0) [dashed] -- (1.35,1.4)  ;
				\draw (0,0) [dashed] -- (1.2,1.8)  ;
				\draw (0,0) [dashed] -- (6.9,1.4)  ;
				\draw (-1, 4) node[above left]{$y$} ;
				\draw (0.9, 1.1) node[above left]{$E_3$} ;
				\draw (0.7, 1) node[below right]{$E_2$} ;
				\draw (6.3, 1.3) node[below]{$E_1$} ;
				\draw (8,0) node[below right]{$g_no$} ;	
				\draw (10,0) node[below right]{$z$} ;	
				\draw (9,4) node[above right]{$g_n x$};
				\draw (4, 4.1) node[above]{$\geq An$};
				\filldraw[black] (8,0) circle (1 pt) ;
				\filldraw[black] (0,0) circle (1 pt) ;
			\end{tikzpicture}
		\end{center}
		\caption{A geometric interpretation of Lemma \ref{geometric lemma}}\label{figure geom interpretation estimate}
	\end{figure}

	The proof of points \ref{geom estimate 1} and \ref{geom estimate 2} is straightforward, so we begin by these. 
	
	\begin{proof}[Proof of estimates \ref{geom estimate 1} and \ref{geom estimate 2}]
		A simple computation gives that 
		\begin{eqnarray}
			(g_n x | g_n o)_o = \frac{1}{2}(b_x (g^{-1}_n o) + d(g_no, o)) \nonumber
		\end{eqnarray}
	Then using assumptions \ref{assum gnx} and \ref{assum drift cat} gives immediately that $(g_nx | g_no )_o \geq (\lambda - \epsilon ) n $, which proves \ref{geom estimate 1}. 	
	
	Now, by definition,  
	\begin{eqnarray}
		(y |g_n o )_o = \frac{1}{2}(d(g_no, o) - b_y (g_no))\nonumber
	\end{eqnarray}
	Then by assumptions \ref{assum y} and \ref{assum drift cat}, we obtain \ref{geom estimate 2}. 
	\end{proof}
	
\begin{figure}
		\centering
		\begin{center}
			\begin{tikzpicture}[scale=1.2]
				\draw (0,0) -- (4,0)  ;
				\draw (2.8,2) -- (2.8,-1)  ;
				\draw (1.2,2) -- (1.2,-1)  ;
				\draw (2,2) -- (2,-1)  ;
				\draw (0,0) node[below left]{$o$} ;
				\draw (4,0) [-stealth] to[bend left = 20] (5, 3);
				\draw (0,0) [-stealth] to[bend right = 20] (-1, 2.5);
				\draw (-1, 2.5) to[bend right = 55] (5,3);
				\draw (-1, 2.5) node[above left]{$y$} ;
				\draw (4,0) node[below right]{$g_no$} ;	
				\draw (5,3) node[above]{$g_n x$};
			\end{tikzpicture}
		\end{center}
		\caption{A "hyperbolic-like" 4 points inequality in Proposition \ref{geometric lemma}.}\label{figure separating hyperplanes}
	\end{figure}
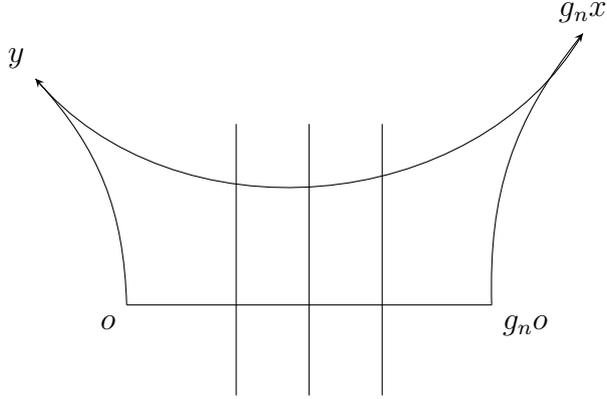

	The proof of point \ref{geom estimate 3} is the hard part. We prove it in two steps. First, we show that under the assumptions, for $n$ large enough, there exist at least three $L$-separated curtains dual to $[o, g_no]$ separating $\{g_n o, g_nx\}$ on the one side and $\{o, y \}$ on the other, see Figure \ref{figure separating hyperplanes}. Then we show that the presence of these curtains implies the result. 
	
	By assumption \ref{assum drift hyp}, for every $n \geq n_0$, $d_L(g_n o, o) \geq An$. For $n\geq n_0$, pick such a $L$-chain separating $o$ and $g_n o $ of size $\geq An$ and define $S(n) \in \mathbb{N}$ as the size of this chain. By Proposition \ref{dual chain}, there exists an $L$-chain dual to $[o, g_no ]$ of size greater than or equal to $\lfloor\frac{S(n)}{4L + 1}\rfloor$ that separates $o$ and $g_no$. Denote by $c_n = \{h^n_i\}_{i=1}^{S'(n)}$ a maximal $L$-chain dual to $[o, g_no]$, separating $o $ and $g_n o $, and orient the half-spaces so that $o \in h_i ^{-}$ for all $i$. When the context is clear, we might omit the dependence in $n$ for ease of notations, and just write $\{h_i\}_{i=1}^{S'(n)}$ for a maximal $L$-chain dual to $[o, g_n o ]$. Recall that $S(n) \geq An$, hence $c_n$ must be of length $S'(n) \geq A'n$, where $A' = \frac{A}{4 L+1}$. 
	
	\begin{lem}\label{estimate y}
		Under the assumptions of Proposition \ref{geometric lemma}, there exists a constant $C$ such that for all $n \in \mathbb{N}$, the number of $L$-separated hyperplanes in $c_n$ that do not separate $\{o, y\}$ and $\{g_n o\} $ is less than $C$. 
	\end{lem}
	\begin{proof}[Proof of Lemma \ref{estimate y}]
		By assumption, $ \{g_n o \}_n $ converges in $(\overline{X_L}, d_L)$ to a point of the boundary $z_L \in \partial X_L$. By Theorem \ref{thm equivariant embedding of boundaries}, there exists an $\iso(X)$-equivariant embedding $\mathcal{I} : \partial X_L \rightarrow \bd X$ that extends the canonical inclusion $X_L \rightarrow X$, and whose image lies in $\mathcal{B}_L$. Denote by $z := \mathcal{I}(z_L)$ the image in $\bd X$ of the limit point $z_L$ by this embedding. 
		
		Denote by $\beta : [0, \infty) \rightarrow X$ a geodesic ray joining $o$ to $z$. Since $z \in \mathcal{B}_L$, there exists an infinite $L$-chain $c = \{k_i\}_{i \in \mathbb{N}}$ that separates $o$ from $z$. Note that because of Lemma \ref{dual chain} and Remark \ref{infinite dual chain}, we can assume that $c$ is a chain of curtains which is dual to the geodesic ray $\beta$. Since $\{g_n o \}_n $ converges in $(\overline{X_L}, d_L)$ to $z_L \in \partial X_L$, and $z$ is the image of $z_L$ by the equivariant embedding $\mathcal{I}$, it implies that $\{g_no\}_n $ converges to $z$ in $X$. The fact that $z \in \mathcal{B}_L$ implies that for all $i \in \mathbb{N}$, there exists $n_0 \in \mathbb{N}$ such that for all $n\geq n_0$, $k_i$ separates $o $ from $g_n o $. Now, we denote by $\gamma : [0, \infty) \rightarrow X$ the geodesic ray that represents $y \in \bd X $. See figure \ref{figure estimate y}. 
		
		\begin{figure}
			\centering
			\begin{center}
				\begin{tikzpicture}[scale=1]
					\draw (0,0) [-stealth]-- (7,0)  ;
					\draw (1.2,4) -- (1.2,-1.5)  ;
					\draw (2.5,4) -- (2.5,-1.5)  ;
					\draw (3.5,4) -- (3.5,-1.5)  ;
					\draw (4.5,4) -- (4.5,-1.5)  ;
					\draw (0,0) node[below left]{$o$} ;
					\draw (2, 3) node[left]{$\gamma$} ;
					\draw (7,0) node[below right]{$z$} ;
					\draw (2, 4.5) node[above left]{$y$} ;
					\draw (2.5, 1) node[above right]{$k_p$} ;
					\draw (3.5, 1) node[above right]{$k_{p+1}$} ;
					\draw (3.5, 0) node[above right]{$r$} ;
					\draw (4.5, 1) node[above right]{$k_{p+2}$} ;
					\draw (5.5,-1) node[below right]{$g_no$} ;
					\draw (6,0) node[above right]{$\beta$} ;
					\filldraw[black] (0,0) circle (1.5pt) ; 
					\draw  (3.3,-0.6) node[below left]{$r'$} ;
					\filldraw[black] (3.5,0) circle (1.5pt) ; 
					\filldraw[black] (3.3,-0.6) circle (1.5pt) ; 
					\filldraw[black] (5.5,-1) circle (1.5pt) ; 
					\draw (0, 0 ) [-stealth] to[bend right]	(2, 4.5);
					\draw (0, 0 ) -- (5.5,-1);
				\end{tikzpicture}
			\end{center}
			\caption{Illustration of Lemma \ref{estimate y}.}\label{figure estimate y}
		\end{figure}
		
		Due to Remark \ref{infinite dual chain}, meeting $c$ infinitely often is equivalent to crossing it, then since $y \neq z$, there exists $p \in \mathbb{N} $ such that $\gamma \subseteq k_p^{-}$. Now consider $n_0$ such that for $n\geq n_0$, $g_n o \in k_{p+2}^{+}$. Fix $n \geq n_0$. Recall that $c_n$ is a maximal $L$-chain dual to $[o, g_n o]$ separating $o$ and $g_n o $. 
		
		Denote by $r \in \beta$ a point in the pole of $k_{p+1}$, and denote by $r'= r'(n)$ the projection of $r$ onto the geodesic $[o, g_n o ]$. Then by Lemma \ref{bottleneck}, 
		\begin{equation*}
			d(o, r'(n)) \leq d(o , r) + 2L +1 \nonumber. 
		\end{equation*} 
		Due to the thickness of the curtains (Remark \ref{remark curtains thick}), the number of curtains in $c_n$ that separate $o $ and $r' (n) $ is $\leq d(o , r) + 2L +2$. We emphasize that this number does not depend on $n\geq n_0$, because for all $n \geq n_0$, $g_n \in k_{p+2}^{+}$ and the previous equation holds.  
		
		Recall that $\gamma \subseteq k_{p}^{-}$, so in particular $\gamma \subseteq k_{p+1}^{-}$. Then by star convexity of the curtains (Lemma \ref{star convexity}), every curtain in $c_n$ whose pole belongs to $[r'(n), g_n o]$ separates $\{o,y\}$ from $g_n o $. Then by the previous argument, the number of curtains that do not separate $\{o, y\}$ from $\{g_no\}$ is less than $d(o , r'(n))$. In particular, the number of curtains that do not separate $\{o, y\}$ from $\{g_no\}$ is less than $d(o, r) + 2L+2 $. Since this quantity does not depend on $n$, we have proven the Lemma. 
	\end{proof}
	
	Now, for a fixed $n$, let us give an estimate for the number of curtains in $ c_n = \{h^n_1 , \dots , h^n_{S'(n)}\} $ that separate $o $ and $g_n x $. When a given $n$ is fixed, we omit the dependence in $n$ and just write $c_n = \{h_1 , \dots , h_{S'(n)}\}$ to ease the notations. Let $\gamma_n : [0, \infty) \rightarrow X$ be the geodesic ray joining $o $ and $g_n x$. Let us take $k_0= k_0(n)$ (depending on $n$) large enough so that for all $k \geq k_0$, 
	\begin{equation}
		| (g_n o | g_n x )_o - (g_no | \gamma_n(k))_o | \leq 1.   \nonumber
	\end{equation}
	
	\begin{lem}\label{estimate gnx}
		Under the assumptions of Proposition \ref{geometric lemma}, the number of $L$-separated hyperplanes in $c_n$ that separate $\{o\}$ and $\{g_n o, \gamma_n(k_0)\} $ is unbounded in $n$. More precisely, for all $M \in \mathbb{N}$, there exists $n_0$ such that for all $n \geq n_0$, the number of $L$-separated hyperplanes in $c_n$ that separate $\{o\}$ and $\{g_n o, \gamma_n(k)\} $ is greater than $M$ for all $k \geq k_0$. 
	\end{lem}
	
	\begin{proof}[Proof of Lemma \ref{estimate gnx}]
		Let $k \geq k_0$. Suppose that the number of curtains in $ c_n = \{h_1 , \dots , h_{S'(n)}\} $ separating $o $ and $\gamma_n(k) $ is less than or equal to $p \in  [0, S'(n)-4]$. Then $\{h_{p+2} , \dots , h_{S'(n)}\} $ is an $L$-chain separating $\{o, \gamma_n(k)\} $ and $\{g_n o\} $. We then denote by $r(n)$ a point on $h_{p+3} \cap [\gamma_n(k), g_no]$ and by $r'(n)$ the projection of $r(n)$ onto $[o, g_n o]$, see Figure \ref{figure estimate gnx}. 
			\begin{figure}[h!]
			\centering
			\begin{center}
				\begin{tikzpicture}[scale=2]
					\draw (0,0) -- (4,0)  ;
					\draw (3.6,2) -- (3.6,-0.5)  ;
					\draw (2.5,2) -- (2.5,-0.5)  ;
					\draw (3,2) -- (3,-0.5)  ;
					\draw (1.5,2) -- (1.5,-0.5)  ;
					\draw (0,0) node[below left]{$o$} ;
					\draw (2.2, 2.5) to[bend right = 40] (4,0);
					\draw (0,0) [-stealth]to[bend right = 20] (2.3,2.8);
					\draw (2.2, 2.5) node[below left]{$\gamma_n(k)$} ;
					\draw (2.3, 2.8) node[above right]{$g_n x$} ;
					\draw (3,0.5) node[above right]{$r(n)$} ;
					\filldraw[black] (2.2, 2.5) circle (1 pt) ;
					\filldraw[black] (3, 0.5) circle (1 pt) ;
					\filldraw[black] (0,0) circle (1 pt) ;
					\filldraw[black] (2.85,0) circle (1 pt) ;
					\filldraw[black] (4,0) circle (1 pt) ;
					\draw (3,-0.5) node[below ]{$h_{p+3}$} ;
					\draw (1.5,-0.5) node[below ]{$h_p$} ;
					\draw (2.85,0) node[above]{$r'(n)$} ;
					\draw (2.5,-0.5) node[below ]{$h_{p+2}$} ;
					\draw (4,0) node[below right]{$g_no$} ;	
				\end{tikzpicture}
			\end{center}
			\caption{Illustration of Lemma \ref{estimate gnx}.}\label{figure estimate gnx}
		\end{figure}
		
		By hypothesis on $k$,
		\begin{eqnarray}
			2((g_nx | g_n o)_o - 1) &\leq& 2(\gamma_n(k)| g_no)_o \nonumber \\
			&=& d(\gamma_n(k) , o ) + d(g_n o , o ) - d(g_n o, \gamma_n (k) ). \nonumber
		\end{eqnarray} 
		Now by the bottleneck Lemma \ref{bottleneck} and the triangular inequality, 
		\begin{eqnarray}
			2(\gamma_n(k)| g_no)_o & = & d(\gamma_n(k) , o ) + d(g_n o , o ) - (d(g_n o, r(n) ) + d(r(n) , \gamma_n(k))) \nonumber \\
			&\leq & d(\gamma_n(k) , o ) + d(g_n o , o ) - (d(g_n o, r' (n)) - (2L + 1) + d(r (n), \gamma_n(k))) \nonumber \\
			&\leq & d(r (n) , o ) + d(g_n o , o ) - d(g_n o, r'(n) ) + 2L + 1 \nonumber \\ 
			&\leq & d(r' (n) , o ) + 2L + 1 + d(g_n o , o ) - d(g_n o, r'(n) ) + 2L + 1 \nonumber \\ 
			& \leq & 2 d(r' (n), o ) + 2(2L + 1). \nonumber
		\end{eqnarray} 
		
		Because the pole of a curtain is of diameter 1, $d(o , r'(n)) \leq d(g_n o, o) - (S'(n) - (p +1))$. However, by assumptions \ref{assum drift hyp} and \ref{assum drift cat} of Proposition \ref{geometric lemma}, one gets that $d(g_no , o ) \leq (\lambda + \varepsilon) n $ and $S(n) \geq An $. Recall that by Lemma \ref{dual chain}, this means that $S'(n) \geq A'n$, where $A' = \frac{A}{4L +10}$. Combining this with the previous result yields
		\begin{eqnarray}
			& & (g_nx | g_n o)_o - 1 \leq d(o, r'(n)) + 2L+1 \nonumber \\
			& \Rightarrow& (\lambda - \varepsilon) n  - 1 \leq (\lambda + \varepsilon )n - (A'n - (p+1)) + 2L+1 \text{ by Lemma \ref{geometric lemma}, \ref{geom estimate 1}}\nonumber  \\
			& \Rightarrow & 0 \leq (2\varepsilon - A' ) n + 2L + p+3.   \nonumber
		\end{eqnarray}
		If $A' > 2 \varepsilon$, there exists $n_0 $ large enough such that for all $n \geq n_0$, the above inequality gives a contradiction. As a consequence, if $A' > 2 \varepsilon$, or equivalently if $ A > 2(4L +10) \varepsilon$, there exists $n_0$ such that for all $n \geq n_0$, the number of curtains in $c_n$ separating ${o } $ and $\{\gamma_n (k), g_no\}$ is greater than $p$.
	\end{proof}

	We can now conclude the proof of Proposition \ref{geometric lemma}. 
	
	\begin{proof}[Proof of estimate \ref{geom estimate 3}]
		Recall that we denote by $\gamma : [0, \infty) \rightarrow X$ the geodesic ray that represents $y \in \bd X $ such that $\gamma(0) = o$ and by $\gamma_n : [0, \infty) \rightarrow X$ the geodesic ray joining $o $ and $g_n x$. Combining Lemma \ref{estimate y} and Lemma \ref{estimate gnx}, we get that if $A > 2(4L+10)\varepsilon$, there exists $n_0$, $k_0$ such that for all $n \geq n_0$ and all $k\geq k_0$, $c_n$ contains at least 3 pairwise $L$-separated curtains that separate $\{o, \gamma(k)\} $ on the one side and $\{g_no, \gamma_n(k)\}$ on the other. Call these hyperplanes $\{h_1, h_2, h_3\}$ and arrange the order so that $h_i \subseteq h_{i+1}^{-}$. Denote by $m_k (n) \in h_2$ some point on the geodesic segment joining $\gamma(k) $ to $\gamma_n(k)$, and $m'_k(n)$ belonging to the geodesic segment $[o, g_no]$ such that $d(m_k(n), m'_k(n)) \leq 2 L + 1$, see Figure \ref{figure proof geom estimate}
		
		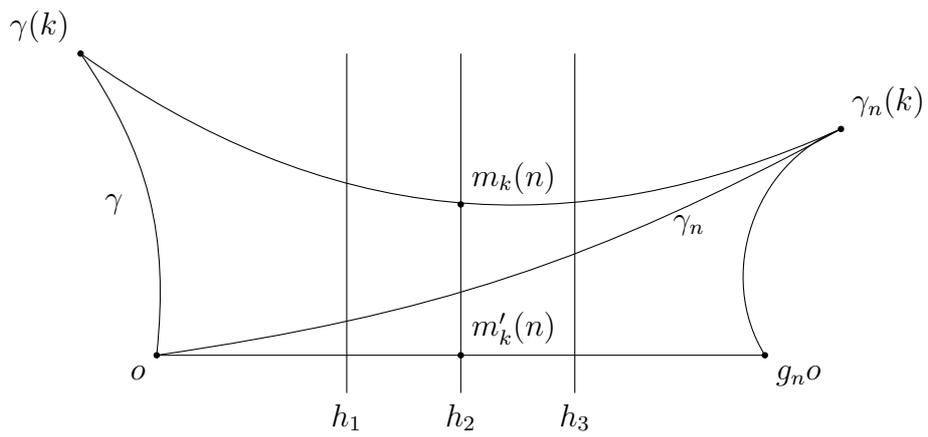
\begin{figure}[h!]
			\centering
			\begin{center}
				\begin{tikzpicture}[scale=1]
					\draw (0,0) -- (8,0)  ;
					\draw (2.5,4) -- (2.5,-0.5)  ;
					\draw (4,4) -- (4,-0.5)  ;
					\draw (5.5, 4) -- (5.5,-0.5)  ;
					\draw (0,0) node[below left]{$o$} ;
					\draw (-0.3, 2) node[left]{$\gamma$} ;
					\draw (0, 0) to[bend right = 20](-1, 4) ;
					\draw (0, 0) to[bend right = 10] (9,3) ;
					\draw (-1, 4) to[bend right = 30] (9,3);
					\draw (8,0) to[bend left = 50] (9,3);
					\draw (-1, 4) node[above left]{$\gamma(k)$} ;
					\draw (2.5, -0.5) node[below]{$h_1$} ;
					\draw (4, -0.5) node[below]{$h_2$} ;
					\draw (5.5, -0.5) node[below]{$h_3$} ;
					\draw (8,0) node[below right]{$g_no$} ;	
					\draw (9,3) node[above right]{$\gamma_n(k)$};
					\draw (4,2) node[above right]{$m_k(n)$};
					\draw (4,0) node[above right]{$m'_k(n)$};
					\draw (7,2) node[below]{$\gamma_n$};
					\filldraw[black] (-1, 4) circle (1 pt) ;
					\filldraw[black] (4,0) circle (1 pt) ;
					\filldraw[black] (4, 2.0) circle (1 pt) ;
					\filldraw[black] (9,3) circle (1 pt) ;
					\filldraw[black] (8,0) circle (1 pt) ;
					\filldraw[black] (0,0) circle (1 pt) ;
				\end{tikzpicture}
			\end{center}
			\caption{Proof of Proposition \ref{geometric lemma}}\label{figure proof geom estimate}
		\end{figure}
		Then we have 
		\begin{eqnarray}
			2(\gamma(k)| \gamma_n(k))_o & = & d(\gamma(k) , o ) + d(o, \gamma_n(k)) - d(\gamma(k), \gamma_n(k)) \nonumber\\
			& \leq & d(\gamma(k) , o ) + d(o, m'_k(n)) + d(m'_k(n), m_k(n) ) \nonumber \\
			& & + \, d(m_k (n), \gamma_n(k)) - d(\gamma(k), \gamma_n(k)) \text{ by the triangular inequality }\nonumber\\
			& \leq & d(\gamma(k) , o ) + d(o, m'_k(n))  - d(\gamma(k), m_k(n)) + 2L + 1 \nonumber
		\end{eqnarray}
		by Lemma \ref{bottleneck}. Since $m'_k(n)$ is on $[o , g_no]$, $d(o , m'_k(n))= d(o, g_no ) - d(g_no , m'_k(n))$. 
		We then have:
		\begin{eqnarray}
			2(\gamma(k)| \gamma_n(k))_o & \leq & d(\gamma(k) , o ) + d(o, g_no ) - d(g_no , m'k(n)) - d(\gamma(k), m_k(n)) + 2L + 1 \nonumber \\
			& = & d(\gamma(k) , o ) + d(o, g_no ) - (d(g_no , m'_k(n)) + d(\gamma(k), m_k(n))) + 2L + 1. \nonumber
		\end{eqnarray}
		 Now observe that 
		 \begin{eqnarray}
		 	d(\gamma(k), g_no) &\leq& d(g_no , m'_k(n)) + d(\gamma(k), m_k(n)) + d(m_k, m'_k(n)) \nonumber \\
		 	& \leq & d(g_no , m'_k(n)) + d(\gamma(k), m_k(n)) + 2L + 1 \text{ by Lemma \ref{bottleneck}}, \nonumber
		 \end{eqnarray}
		 hence $d(\gamma(k), g_n o ) - (2 L + 1) \leq d(g_no , m'_k(n)) + d(\gamma(k), m_k(n))$. Then 
		\begin{eqnarray}
			2(\gamma(k)| \gamma_n(k))_o & \leq & d(\gamma(k) , o ) + d(o, g_no ) - (d(\gamma(k), g_n o ) - (2 L + 1)) + 2L + 1 \nonumber \\
			& = & d(\gamma(k) , o ) + d(o, g_no ) - d(\gamma(k), g_n o ) + 2(2 L + 1) \nonumber \\
			& = & 2(\gamma(k)| g_n o )_o + 2(2 L + 1) \nonumber.
		\end{eqnarray}
		As $k \rightarrow \infty$, one obtains that $(g_nx| y)_o \leq (g_no |y)_o + (2 L + 1)$, and the result follows from \ref{estimate 2}. 
	\end{proof}

	\subsection{Proof of the Central Limit Theorem}
	
	In this section, we prove the main result of the paper. Let $G$ be a discrete group and $G \curvearrowright X$ a non-elementary action by isometries on a proper $\cat$(0) space $X$. Let $\mu \in \prob(G) $ be an admissible probability measure on $G$ with finite second moment, and assume that $G $ contains a rank one element. Let $o \in X$ be a basepoint of the random walk. Let $\lambda$ be the (positive) drift of the random walk provided by Theorem \ref{drift thm}.
	We assume the action on $X$ to be non elementary and rank one, hence due to Proposition \ref{prop non elem}, there exists a number $L \geq 0$ such that $G$ acts by isometries on $X_L = (X, d_L) $ non elementarily. Then one can consider the random walk $(Z_n(\omega) o)_n$ as a random walk on $(X, d_L) $, which we will write $(Z_n \tilde{o})_n $ when the context is not clear. The model $(X, d_L) $ is hyperbolic, so we can apply the results of Maher and Tiozzo \cite{maher_tiozzo18} summarized in Section \ref{results hyp}. In particular, due to Theorem \ref{mt conv} (along with Gou\"el's result recalled here in Remark \ref{rem gouezel cv non-sep} for non-separable hyperbolic spaces), the random walk $(Z_n \tilde{o})_n$ in $X_L$ converges to a point of the Gromov boundary $\partial X_L$ of $(X, d_L)$. 
	
	Moreover, since we assume $\mu $ to have finite first moment (for the action on the $\cat$(0) space $X$), and since $d(x, y ) +1 \geq d_L(x, y )$ for all $x, y \in X$, the measure $\mu $ is also of finite first moment for the action on the hyperbolic model $(X, d_L)$. In particular, the drift $\tilde{\lambda}$ of the random walk $(Z_n \tilde{o})_n$ is almost surely positive. In other words, we have that $\mathbb{P}$-almost surely, 
	\begin{equation}
		\lim_{n \rightarrow \infty} \frac{1}{n} d(Z_n(\omega) \tilde{o}, \tilde{o}) = \tilde{\lambda} >0. \nonumber
	\end{equation}
	
	Due to Theorem \ref{measure thm cat0}, there exists a unique $\mu$-stationary probability measure $\nu$ on $\overline{X}$. If we define $\mui \in \prob(G)$ by $\mui (g) = \mu(g^{-1})$, $\mui $ is still admissible and of finite second moment. We denote by $\nui$ the unique $\mui$-stationary measure on $\overline{X}$. 
	
	We recall that the Busemann cocycle $\beta : G \times \overline{X} \rightarrow \mathbb{R}$ is defined by:
	\begin{equation*}
		\beta(g,x) = b_x (g^{-1} o).
	\end{equation*}
	
	Our goal is to apply Theorem \ref{criterium clt benoist quint} to the Busemann cocycle $\beta$. The results of Section \ref{strategy} show that proving a central limit theorem for the random walk $(Z_n(\omega)o)_n$ amounts to proving that $\beta$ is centerable. As in the works of \cite{benoist_quint16}, \cite{horbez18} and \cite{fernos_lecureux_matheus21}, the natural candidate to solving the cohomological equation \eqref{cohom equation} is the function:
	\begin{equation*}
		\psi ( x ) = -2 \int_{\overline{X}} (x | y)_o d\nui(y).
	\end{equation*}
	
	\begin{prop}\label{prop centerable cocycle}
		Let $G$ be a discrete group and $G \curvearrowright X$ a non-elementary action by isometries on a proper $\cat$(0) space $X$. Let $\mu \in \prob(G) $ be an admissible probability measure on $G$ with finite second moment, and assume that $G $ contains a rank one element. Let $o \in X$ be a basepoint of the random walk. Then the Borel map
		$\psi(x) = \int_{\overline{X}} (x |y)_o \, d\nui(y)$ is well-defined and essentially bounded. 
	\end{prop}

	In order to show that $\psi$ is well-defined and bounded, we need the following statement, which resembles \cite[Proposition 4.2]{benoist_quint16}.

	\begin{prop}\label{estimates}
		Let $G$ be a discrete group and $G \curvearrowright X$ a non-elementary action by isometries on a proper $\cat$(0) space $X$. Let $\mu \in \prob(G) $ be an admissible probability measure on $G$ with finite second moment, and assume that $G $ contains a rank one element. Let $o \in X$ be a basepoint for the random walk $(Z_n(\omega)o)_n$. Let $\lambda$ be the (positive) drift of the random walk, and $\nu$ a $\mu$-stationary measure on $\overline{X}$. Assume that there exists $a >0 $ and $(C_n)_n \in \ell^1 (\mathbb{N})$ such that for almost every $x, y \in \overline{X}$, we have, for every $n$: 
		\begin{enumerate}
			\item $\mathbb{P}((Z_n o | Z_n x )_o \leq an) \leq C_n $; \label{estimate 1}
			\item $\mathbb{P}((Z_n o | y )_o \geq an) \leq C_n $; \label{estimate 2}
			\item $\mathbb{P}((Z_n x | y )_o \geq an) \leq C_n $. \label{estimate 3}
		\end{enumerate}
		Then one has: 
		\begin{eqnarray}
			\sup_{x \in \overline{X}} \int_{\overline{X}} (x |y)_o d\nu(y) < \infty \nonumber.
		\end{eqnarray}
	\end{prop}
	\begin{proof}
		Suppose that there exist $a>0$, $(C_n)_n \in \ell^1 (\mathbb{N})$ such that for almost every $x, y \in \overline{X}$, we have estimates \ref{estimate 1}, \ref{estimate 2} and \ref{estimate 3}. We get: 
		\begin{eqnarray}
			\nu(\{x \in X | (x|y) \geq an \}) & = & \int_{\overline{X}} \mu^{\ast n }(\{g\in G \, | \, (gx |y )_o \geq an \})d\nu(x) \text{ by $\mu$-stationarity} \nonumber \\
			& \leq & \int_{\overline{X}} C_n d\nu(x) = C_n \text{ by estimate \ref{estimate 3}}.\nonumber 
		\end{eqnarray}
		Then, define $A_{n, y } := \{ x \in \overline{X} \, | \, (x|y )_o \geq an\}$, so that by splitting along the subsets $A_{n-1, y } - A_{n, y }$, one gets 
		\begin{eqnarray}
			\int_{\overline{X}} (x |y)_o d\nu(x) & \leq & \sum_{n \geq 1} a n (\nu(A_{n-1, y }) - \nu(A_{n, y })) \nonumber \\
			& \leq &  \sum_{n \geq 1} an (C_{n-1} - C_n) \nonumber \\
			& = & a + \sum_{n \geq 1} aC_n (n+1 - n) < \infty. \nonumber
		\end{eqnarray}
	\end{proof}
	
	We want to show that estimates from Proposition \ref{estimates} hold. As we will see, estimates \ref{estimate 1} and \ref{estimate 2} are quite straightforward to check using the positivity of the drift. Most of the work concerns estimate \ref{estimate 3}.

	Combining Proposition \ref{Busemann cocyle prop} with Theorem \ref{drift thm} and \cite[Proposition 3.2]{benoist_quint16CLTlineargroups}, one obtains the following:
	
	\begin{prop}\label{proba control}
		Let $G$ be a discrete group and $G \curvearrowright X$ a non-elementary action by isometries on a proper $\cat$(0) space $X$. Let $\mu \in \prob(G) $ be an admissible probability measure on $G$ with finite second moment, and assume that $G $ contains a rank one element. Let $o \in X$ be a basepoint of the random walk. Let $\lambda$ be the (positive) drift of the random walk. Then, for every $\varepsilon >0 $, there exists $(C_n)_n \in \ell^1 (\mathbb{N})$ such that for any $x \in \overline{X}$, 
		\begin{eqnarray}
			& & \mathbb{P}(|\beta(Z_n, x) - n \lambda | \geq \varepsilon n) \leq C_n;\label{prob estimate 1}\\
			& & \mathbb{P}(|\beta(Z^{-1}_n, x) - n \lambda | \geq \varepsilon n) \leq C_n ;  \label{prob estimate 2} \\
			& & \mathbb{P}(|d(Z_n o, o) - n \lambda | \geq \varepsilon n) \leq C_n.  \label{prob estimate 3}
		\end{eqnarray}
	\end{prop}
		
	\begin{proof}
		Recall that by Proposition \ref{Busemann cocyle prop}, $\beta $ is a continuous cocycle such that
		\begin{eqnarray}
			& & \int_G \sup_{x \in \overline{X}} |\beta(g,x)|^2 d\mu(g) < \infty \text{ and }
			\int_G \sup_{x \in \overline{X}} |\beta(g,x)|^2 d\mui(g) < \infty. \nonumber
		\end{eqnarray}
		Moreover, 
		\begin{eqnarray}
			\lambda = \int_{G \times \overline{X}} \beta(g, x) d\mu(g)d\nu(x)= \int_{G \times \overline{X}} \beta(g, x) d\mui(g)d\nui(x). \nonumber
		\end{eqnarray} 
		We can then apply \cite[Proposition 3.2]{benoist_quint16CLTlineargroups}: for every $\varepsilon >0 $, there exists a sequence $(C_n) \in \ell^1 (\mathbb{N})$ such that for every $x \in \overline{X}$, 
		\begin{eqnarray}
			\mathbb{P}\big(\omega \in \Omega \, : \, \big|\frac{\beta(Z_n(\omega), x)}{n} - \lambda\big| \geq \epsilon\big) \leq C_n \nonumber. 
		\end{eqnarray}
		The same goes for $\mui $ and $\nui$, which gives estimates \eqref{prob estimate 1} and \eqref{prob estimate 2}. 
		
		Estimate \eqref{prob estimate 3} is then a straightforward consequence of Proposition \ref{approx displacement horo}. 
	\end{proof}

	The following Lemma will also be important in the proof of Proposition \ref{prop centerable cocycle}. 
	
	\begin{lem}\label{Lemma drift hyp}
		Let $G$ be a discrete group and $G \curvearrowright X$ a non-elementary action by isometries on a proper $\cat$(0) space $X$. Let $\mu \in \prob(G) $ be an admissible probability measure on $G$ with finite second moment, and assume that $G $ contains a rank one element. Let $o \in X$ be a basepoint of the random walk. Let $\lambda$ be the (positive) drift of the random walk.
		
		Then there exist $L >0$, $\lambda_L > 0$ such that almost surely, $\liminf_n \frac{d_L (Z_n o, o)}{n} = \lambda_L$. Moreover, there exists $A > 0$ and $(C_n) \in \ell^1 (\mathbb{N}) $ such that 
		\begin{eqnarray}
			\mathbb{P}\big( d_L (Z_n o, o) < An\big) \leq C_n \nonumber. 
		\end{eqnarray}
	\end{lem}

	\begin{proof}
		The action $G \curvearrowright (X,d)$ is non elementary and contains a rank one element, hence by Proposition \ref{prop non elem}, there exists $L$ such that the action $G \curvearrowright (X, d_L)$ is non-elementary as the loxodromic isometries $g$ and $h$ are independent. We can then apply Theorem \ref{mt drift}, which gives the Lemma. 
	\end{proof}

	Let us now complete the proof of Proposition \ref{prop centerable cocycle}.
	\begin{proof}[Proof of Proposition \ref{prop centerable cocycle}]
		By assumptions, we can apply Theorem \ref{mt conv}: there exists $L >0 $ such that $(Z_n(\omega) o )_n$ converges in $(X_L, d_L)$ to a point $z_L$ of the boundary. By Theorem \ref{measure thm cat0}, there is a unique $\mu$-stationary measure $ \nu$ on $\bd X$, and this measure is non-atomic.
		 
		Fix $A$ as in Lemma \ref{Lemma drift hyp}, and $(C_n)_n \in \ell^1 (\mathbb{N})$ such that
		\begin{eqnarray}
			\mathbb{P}\big( d_L (Z_n o, o) < An\big) < C_n \nonumber. 
		\end{eqnarray} 
		Now take $0 < \varepsilon < \min(\frac{A}{2(4L+10)}, \lambda/2)$. Due to Proposition \ref{proba control}, there exists a sequence $C'_n \in \ell^1( \mathbb{N})$ such that 
		\begin{eqnarray}
			& & \mathbb{P}(|\beta(Z_n, x) - n \lambda | \geq \varepsilon n) \leq C'_n\nonumber\\
			& & \mathbb{P}(|\beta(Z^{-1}_n, x) - n \lambda | \geq \varepsilon n) \leq C'_n  \nonumber \\
			& & \mathbb{P}(|d(Z_n o, o) - n \lambda | \geq \varepsilon n) \leq C'_n. \nonumber 
		\end{eqnarray}
		We can assume that $C_n = C'_n $ for all $n$. Then for $\nu$-almost every $x, y \in \bd X$, we have the quantitative assumptions in Proposition \ref{geometric lemma}: 
		\begin{enumerate}[label= (\roman*)]
			\item $ \{Z_no  \}_n $ converges in $(\overline{X_L}, d_L)$ to a point of the boundary $z_L \in \partial X_L$, whose image in $\bd X $ by the embedding $\mathcal{I}$ is not $y$; 
			\item $\mathbb{P}\big( d_L(Z_n o, o) \geq An \big) \geq 1 - C_n$; 
			\item $\mathbb{P}\big(|b_x (Z_n^{-1} o ) - n \lambda | \leq \varepsilon n\big) \geq 1-C_n$; 
			\item $\mathbb{P}\big(|b_y (Z_n o ) - n \lambda | \leq \varepsilon n\big) \geq 1- C_n$; 
			\item $\mathbb{P}\big(|d(gZ_n o , o) - n \lambda | \leq \varepsilon n\big) \geq 1- C_n$. 
		\end{enumerate}	
		As a consequence, one obtains that for $\nu$-almost every $x, y \in \bd X$, the probability that these estimates are not satisfied is bounded above by $4 C_n$. Now choosing $a \in (\varepsilon, \lambda - \varepsilon)$, we get that for $n$ large enough,
		\begin{enumerate}
			\item $\mathbb{P}\big((g_nx |g_n o )_o \geq a n\big) \geq 1 - 4C_n$; 
			\item $\mathbb{P}\big((y | g_n o )_o \leq a n\big)\geq 1 - 4C_n$;
			\item $\mathbb{P}\big((y| g_n x )_o \leq a n\big) \geq 1 - 4 C_n$.
		\end{enumerate}
	
	Since the sequence $(4C_n)_n$ is still summable, we can apply Proposition \ref{estimates}, that states that the function $\psi $ defined by 
	\begin{equation*}
		\psi ( x ) = -2 \int_{\overline{X}} (x | y)_o d\nui(y)
	\end{equation*}
	is well-defined, and Borel by Fubini. Moreover, $\psi $ is essentially bounded:
	\begin{eqnarray}
		\sup_{x \in \overline{X}} \int_{\overline{X}} (x |y)_o d\nu(y) < \infty \nonumber.
	\end{eqnarray} 
	\end{proof}
	
	\begin{cor}\label{cor centerable cocycle}
		Under the same assumptions as in Proposition \ref{prop centerable cocycle}, the cocycle $\beta(g, x) = b_x (g^{-1} o )$ is centerable.
	\end{cor}
	
	\begin{proof}
		By Proposition \ref{prop centerable cocycle}, the function $\psi $ defined by 
		\begin{equation*}
			\psi ( x ) = -2 \int_{\overline{X}} (x | y)_o d\nui(y).
		\end{equation*}
		is well-defined, Borel and essentially bounded. Also, as observed in \cite[Lemma 1.2]{benoist_quint16}, a quick computation shows that for all $g \in G$, $x,y \in \overline{X}$:
		\begin{eqnarray}
			b_x(g^{-1}o) = -2(x| g^{-1}y)_o + 2(gx|y)_o + b_y (go) \nonumber. 
		\end{eqnarray}
		Fix $x \in X$. Integrate this equality on $(G \times \bd X, \mu \otimes \nui)$ gives 
		\begin{eqnarray}
			\int_G  \beta(g,x) d\mu(g) &=& -2\int_G \int_{\bd X}(x| g^{-1}y)_o d\mu(g) d\nui(y) \nonumber \\
			& & + 2\int_G \int_{\bd X}(gx|y)_o d\mu(g) d\nui(y)+ \int_G \int_{\bd X} \beta(g^{-1}, x)d\mu(g) d\nui(y) \nonumber \\
			&  = & -2\int_G \int_{\bd X}(x| gy)_o d\mui(g) d\nui(y) - \int_G \psi(gx)d\mu(g) \nonumber \\
			& & +\int_G \int_{\bd X} \beta(g, x)d\mui(g) d\nui(y) \nonumber.
		\end{eqnarray}
		But $\int_G \int_{\bd X}(x| gy)_o d\mui(g) d\nui(y)=  \int_{\bd X}(x| y)_o d\nui(y)$ because $\nui$ is $\mui$-stationary. Also, by point \ref{average Busemann} in Proposition \ref{Busemann cocyle prop}, we have that 
		\begin{eqnarray}
			\int_G \int_{\bd X} \beta(g, x)d\mui(g) d\nui(y) = \lambda.
		\end{eqnarray}
		Combining these, we get: 
		\begin{eqnarray}
			\int_G \beta(g,x) d\mu(g) = \psi(x) - \int_G \psi(gx)d\mu(g) + \lambda. \nonumber
		\end{eqnarray}
		
		Hence if we define $\beta_0 (g, x) = \beta(g, x) - \psi(x) + \psi(gx) $, we obtain that for all $x \in \overline{X}$, 
		\begin{eqnarray}
			\int_G \beta_0(g, x) d\mu(g)= \lambda, 
		\end{eqnarray}
		and the cocycle $\beta_0$ has constant drift $\lambda$. Then by Remark \ref{rem average cocycle}, $\beta$ is centerable with average $\lambda$, as wanted. 
	\end{proof}
	
	We can now state the following. 
	\begin{thm}\label{thm cv gaussian law}
		Let $G$ be a discrete group and $G \curvearrowright X$ a non-elementary action by isometries on a proper $\cat$(0) space $X$. Let $\mu \in \prob(G) $ be an admissible probability measure on $G$ with finite second moment, and assume that $G $ contains a rank one element. Let $o \in X$ be a basepoint of the random walk. Let $\lambda$ be the (positive) drift of the random walk. Then the random variables $(\frac{1}{\sqrt{n}}(d(Z_n o, o) - n \lambda))_n $ converge in law to a Gaussian distribution $N_\mu$. Furthermore, the variance of $N_\mu$ is given by 
		\begin{eqnarray}
			\int_{G \times \bd X} (b_x(g^{-1} o) - \psi(x) + \psi(gx) - \lambda)^2 d\mu(g)d\nu(x).
		\end{eqnarray}
	\end{thm}

	\begin{proof}
		By Corollary \ref{cor centerable cocycle}, the cocycle $\beta$ is centerable, with average $\lambda$. Since the measure $\nu$ is the unique $\mu$-stationary measure on $\overline{X}$, we can then apply Theorem \ref{criterium clt benoist quint}: the random variables $(\frac{1}{\sqrt{n}}(\beta(Z_n(\omega),x) - n \lambda))_n$ converge to a Gaussian law $N_\mu $. But thanks to Proposition \ref{approx displacement horo}, this is equivalent to the convergence of the random variables $(\frac{1}{\sqrt{n}}d(Z_n(\omega)o, o) - n \lambda)_n$ to a Gaussian law. Moreover, by Theorem \ref{criterium clt benoist quint} and Proposition \ref{Busemann cocyle prop}, the covariance 2-tensor of the limit law is given by 
		
		\begin{eqnarray}
			\int_{G \times \bd X} (\beta_0(g, z) - \lambda)^2 d\mu(g) d\nu(z) \nonumber, 
		\end{eqnarray}
		where $\beta_0 (g,x)= \beta(g,x) - \psi(x) + \psi(gx)$. This yields the result. 
	\end{proof}

	In order to prove Theorem \ref{thm clt cat0}, it only remains to prove that the limit law is non-degenerate. This is what we do in the next Proposition. 
	
	\begin{prop}
		With the same assumptions and notations as in Theorem \ref{thm cv gaussian law}, the covariance 2-tensor of the limit law satisfies:
		\begin{eqnarray}
			\int_{G \times \bd X} (\beta_0(g, z) - \lambda)^2 d\mu(g) d\nu(z) >0 \nonumber. 
		\end{eqnarray}
	In particular, the limit law $N_\mu$ of the random variables $(\frac{1}{\sqrt{n}}(d(Z_n o, o) - n \lambda))_n $ is non-degenerate. 
	\end{prop}
	
	\begin{prop}\label{prop non deg limit law cat}
		With the same assumptions and notations as in Theorem \ref{thm cv gaussian law}, the covariance 2-tensor of the limit law satisfies:
		\begin{eqnarray}
			\int_{G \times \bd X} (\beta_0(g, z) - \lambda)^2 d\mu(g) d\nu(z) >0 \nonumber. 
		\end{eqnarray}
		In particular, the limit law $N_\mu$ of the random variables $(\frac{1}{\sqrt{n}}(d(Z_n o, o) - n \lambda))_n $ is non-degenerate. 
	\end{prop}
	
	In the course of the proof, we shall use the following fact. We give the proof for completeness. 
	
	\begin{lem}\label{lem attracting pt supp}
		We use the same assumptions and notations as in Theorem \ref{thm cv gaussian law}. Let $g \in G$ be a contracting isometry of $X$, and let $\xi^+ \in \bd X$ be its attracting fixed point at infinity. Then $\xi^+ \in \supp(\nu)$, where $\nu$ is the unique $\mu$-stationary measure on $\overline{X}$. 
	\end{lem}
	\begin{proof}
		Denote by $\xi^-  \in \bd X$ the repelling fixed point in of $g$. The isometry $g$ is contracting, hence by \cite[Lemma 4.4]{hamenstadt09} it acts on $\bd X$ with North-South dynamics. This means that for every neighbourhood $U$ of $\xi^{+}$, $V$ of $\xi^{-}$ in $\bd X$, there exists $k$ such that for all $n \geq k$, $g^n (\bd X  - V) \subseteq U $ and $g^{-n} (\bd X  - U) \subseteq V $. It is standard that $\nu$ is non-atomic, see for instance \cite[Lemme 4.2.6]{le-bars23these}. Hence there exists a neighbourhood $V$ of $\xi ^{-}$ such that $\nu (\bd X - V) > 0$. Fix such a $V$, and let $U$ be any neighbourhood of $\xi^{+}$. Take $k$ large enough so that for all $n \geq k$, $g^n (\bd X  - V) \subseteq U $. Since $\mu$ is admissible, there exists $p' \in \mathbb{N}$ such that $g^k \in \supp(\mu^{\ast p'})$. Check that $\nu$ is still $\mu^{\ast p'}$-stationary, therefore 
		\begin{eqnarray}
			\sum_{h \in G} \nu(h^{-1}U) \mu^{\ast p'} (h) = \nu(U) \nonumber. 
		\end{eqnarray}
		In particular, by North-South dynamics, 
		\begin{eqnarray}
			\nu (U) \geq \nu (g ^{-k}U) \mu^{\ast p'} (g^k) \geq \nu (\bd X - V) \mu^{\ast p'} (g^k) > 0 \nonumber. 
		\end{eqnarray} 
		This is true for every neighbourhood $U$ of $\xi^{+}$, hence $\xi^{+} \in \supp(\nu)$. 
	\end{proof}

	\begin{proof}[Proof of Proposition \ref{prop non deg limit law cat}]
		Let $g$ be a contracting isometry in $G$. Recall that $g$ has an axis $\gamma \subseteq X$ on which $g$ acts as a translation, and let $\xi^{+}, \xi^{-}$ be its attracting and repelling fixed points in $\bd X$ respectively. We let $l(g) = \lim \frac{d(g^n o,o)}{n}$ be the translation length of $g$ in $(X,d)$. Observe that 
		\begin{eqnarray}
			l(g) = \lim_n \frac{b_{\xi^{+}}(g^{-n} o)}{n}, \label{translation cocycle}
		\end{eqnarray}
		where $b_{\xi^{+}} $ is the horofunction centered on $\xi^+$ and based at $o $. 
		Indeed, if $o$ belongs to $\gamma$, then $b_{\xi^{+}}(g^{-n} o) = d(g^n o , o)$ and equation \eqref{translation cocycle} is true. If $o$ does not belong to $\gamma$, take $o' \in \gamma$, and by triangular inequality, 
		\begin{eqnarray}
			|b_{\xi^{+}}(g^{-1} o) - b_{\xi^{+}}^{o'}(g^{-1}o') |\leq 2 d(o, o') \nonumber,
		\end{eqnarray}
		where $b_{\xi^{+}}^{o'}$ is the horofunction with basepoint $o'$. Since $l(g) = \lim_n \frac{1}{n} b_{\xi^{+}}^{o'}(g^{-n} o')$, we obtain that $l(g) = \lim_n \frac{1}{n} b_{\xi^{+}}(g^{-n} o)$. 
		\newline 
		
		Suppose by contradiction that $\int_{G \times \bd X} (\beta_0(h, z) - \lambda)^2 d\mu(h) d\nu(z) = 0$. This means that for almost every $\xi \in \supp(\nu)$ and $h \in \supp(\mu)$, 
		\begin{eqnarray}
			b_\xi (h^{-1} o) - \lambda = \psi (\xi) - \psi(h\xi) \nonumber. 
		\end{eqnarray} 
		Since $\psi $ is bounded and continuous, we get that for every $\xi \in \supp(\nu)$ and every $h \in \supp(\mu)$, $|b_\xi (h^{-1} o) - \lambda| \leq 2 \| \psi \|$. 
		\newline 
		
		Now consider the random walk generated by $\mu^{\ast p}$, for $p \geq 1$. Observe that $\mu^{\ast p }$ is still admissible of finite second moment and that $\nu$ is still a $\mu^{\ast p }$-stationary measure on $\bd X$. We can then apply Theorems \ref{drift thm} and \ref{thm cv gaussian law}, so that the random walk generated by $\mu^{\ast p}$ converges to the boundary with positive drift $l_X (\mu^{\ast p }) = p \lambda >0$. By the previous argument, for almost every $\xi \in \supp(\nu)$ and every $h \in \supp(\mu^{\ast p})$, 
		\begin{eqnarray}\label{eq non deg}
			|b_\xi (h^{-1} o) - p\lambda| \leq 2 \| \psi \|. 
		\end{eqnarray} 
		
		Let $g$ be a contracting element in $G$, and let $\xi^+ $ be its attracting fixed point. Because $\mu $ is admissible, there exists $m $ such that $\mu^{\ast m} ( g) > 0$. Then by Equation \eqref{eq non deg}, for all $n \geq 1$, $|b_{\xi^{+}}(g^{-n}o) - nm \lambda | \leq 2 \| \psi \|$. By Lemma \ref{lem attracting pt supp}, we can apply equation \eqref{translation cocycle}, and we obtain that 
		\begin{eqnarray}
			\lim_n \frac{b_{\xi^{+}}(g^{-n}o)}{n} =  l(g) = m \lambda \nonumber. 
		\end{eqnarray}
		But there also exists $q \in \mathbb{N}^\ast$ such that $1 \in \supp(\mu^{\ast q})$, hence $g \in \supp(\mu^{\ast (m+q)})$ and by the same argument, $l(g) = (m + q) \lambda$. Since by Theorem \ref{drift thm}, $\lambda$ is positive, we get a contradiction. 
	\end{proof}

	\bibliographystyle{alpha}
	\bibliography{bibliography}
\end{document}